\documentclass[a4wide,11pt]{article}
\usepackage{geometry}                

\usepackage{graphicx}
\usepackage{authblk}
\usepackage{amssymb}
\usepackage{epstopdf}
\usepackage[sans]{dsfont}
\usepackage[T1]{fontenc}
\usepackage[english]{babel}
\usepackage{latexsym}
\usepackage{mathrsfs}
\usepackage{amscd}
\usepackage{graphicx}
\usepackage{color}
\usepackage{float}
\usepackage{amsmath}
\usepackage{amsfonts}
\numberwithin{equation}{section}
	\usepackage{enumerate}
\usepackage{amsthm}
\usepackage[numbers]{natbib}
\usepackage[bookmarks=true,colorlinks=true,linkcolor={blue},urlcolor={blue}, citecolor={blue},pdfstartview={XYZ null null 1.22}]{hyperref}%

\def\be#1\ee{\begin{equation}#1\end{equation}}

\newtheorem{thm}{Theorem}[section]

\newtheorem{prop}[thm]{Proposition} 

\newtheorem{rem}{Remark} 

\newcommand{\bq}{\begin{equation}}
\newcommand{\eq}{\end{equation}}

\def\bqa{\begin{eqnarray}}
\def\eqa{\end{eqnarray}}

\def\tw{\tilde w}

\newcommand{\bd}{\begin{displaymath}}
\newcommand{\ed}{\end{displaymath}}
\newcommand{\ba}{\begin{eqnarray}}
\newcommand{\ea}{\end{eqnarray}}

\begin{document}

\title{Boltzmann type control of opinion consensus through leaders}

\author[1]{G. Albi \thanks{giacomo.albi@unife.it}}
\author[1]{L. Pareschi \thanks{lorenzo.pareschi@unife.it}}
\author[1]{M. Zanella \thanks{mattia.zanella@unife.it}}
\affil[1]{University of Ferrara, Department of Mathematics and Computer Science, Via Machiavelli 35, I-44121 Ferrara, ITALY}

\maketitle

%
%
%
%
%

\begin{abstract}
The study of formations and dynamics of opinions leading to the so called opinion consensus is one of the most important areas in mathematical modeling of social sciences. Following the Boltzmann type control recently introduced in \cite{AHP13}, we consider a group of opinion leaders which modify their strategy accordingly to an objective functional with the aim to achieve opinion consensus. 
The main feature of the Boltzmann type control is that, thanks to an instantaneous binary control formulation, it permits to embed the minimization of the cost functional into the microscopic leaders interactions of the corresponding Boltzmann equation. 
The related Fokker-Planck asymptotic limits are also derived which allow to give explicit expressions of stationary solutions.
The results demonstrate the validity of the
Boltzmann type control approach and the capability of the leaders control to strategically lead the followers opinion.
\end{abstract}

\section{Introduction}
Mean-field games and mean-field type control theory has raised a lot of interest in the recent years (see for example \cite{BFY:13,BDMW14,DLR14,FS13,LL07} and the references therein). The general setting consists in a control problem involving a very large number of agents where both the evolution of the state and the objective functional of each agent are influenced by the collective behaviour of all other agents. Typical examples in socio-economical sciences, biology and engineering are represented by the problems of persuading voters to vote for a specific candidate, influencing buyers towards a given good or asset, forcing human crowds or group of animals to follow a specific path or to reach a desired zone, or optimizing the traffic flow in road networks and supply chains\cite{AlbiPareschi2013aa,Ben:05,BDMW14,Borzi14,CPT11,DMPW:09,HertyRinghofer2011ac,HLL05,LuMa:99,MaPa:12}.

In this paper we focus on control problems where the collective behaviour corresponds to the formation process of opinion consensus\cite{Ben03,BS:09,GGS:82,MoTa:13,Sznajd00,Tos06,HK02,NPT:10}. In particular, we consider models where the control strategy is based on hierarchical leadership. This hierarchical leadership
concept is discussed in \cite{DMPW:09}, where a population of leaders is considered giving rise to aggregate opinions and
convergence towards specific patterns. Opinion dynamics in presence of different populations has been previously introduced in \cite{BA07,VKR03,Weisbuch05}. We mention here that control through leaders in self-organized flocking systems has been studied in \cite{AlbiPareschi2013aa,Borzi14}. 

We introduce a hierarchical opinion formation dynamics where the leaders aim at controlling the followers through a suitable cost function which characterizes the leaders strategy in trying to influence the followers opinion. Based on
this microscopic model, we develop a Boltzmann type optimal control approach following the ideas recently presented in \cite{AHP13}. The approach is closely related to model predictive and instantaneous control techniques\cite{MayneMichalska1990aa, CaBo:04,HT83,CHK99}. We derive an explicit controller for the leader dynamic using an instantaneous binary control framework on the microscopic level
and, similarly to the mean field control, study the related kinetic description for large number of agents. Thanks to this formulation, the minimization of the cost functional is embedded into the microscopic leaders interactions of the corresponding Boltzmann equation. 


The rest of the manuscript is organized as follows. First in Section \ref{sec:Binary} we introduce the microscopic model of the leader strategy in the leader-follower interactions and derive the corresponding Boltzmann type control formulation. The main properties of the kinetic model are studied in Section \ref{sec:Boltzmann}, in particular we show that  that the leaders control strategy may lead the followers opinion towards the desired state. Explicit asymptotic opinion distributions are computed in Section \ref{sec:FP} using an approximated Fokker-Planck description derived in the so-called quasi invariant opinion limit. Several numerical results confirm the theoretical analysis in Section \ref{sec:Numeric}.  

\section{Microscopic models of opinion control through leaders}\label{sec:Binary}
A rather common assumption in opinion formations is that interactions are formed mainly by binary exchange of informations, see for example \cite{BS:09,GGL:12,PT:13,Tos06}. Similar to \cite{DMPW:09} we are interested in the opinion formation process of a followers' population steered by the action of a leaders' group. The major novelty here is that the leaders' behaviour is driven by a suitable control strategy based on the interplay between the desire to force followers towards a given state and the necessity to keep a position close the the mean opinion of the followers in order to influence them.
In the following  we first generalize the approach of \cite{AHP13} starting from a differential system describing the evolution of the two populations of leaders and followers.
 In the second part we present a binary interaction model for the same dynamic showing how the two descriptions are related.

\subsection*{Microscopic modeling }
We assume to have two populations, one of followers and one of leaders. Each follower is mutually influenced by the other followers and by the leaders, whose target is to steer the followers' opinion to a desired configuration of consensus following some prescribed strategy. We consider the evolution of a population of $N_L$ leaders and $N_F$ followers, with opinions $w_i, \tilde{w}_k\in I=[-1,1]$, for $i=1,\ldots,N_F$ and $k=1,\ldots,N_L$, evolving according to
\begin{align}
&\dot{w}_{i}=\dfrac{1}{N_F}\displaystyle\sum_{j=1}^{N_F}P\left(w_{i},w_{j}\right)\left(w_{j}-w_{i}\right)+\dfrac{1}{N_{L}}\sum_{h=1}^{N_{L}}S\left(w_{i},\tw_{h}\right)\left(\tw_{h}-w_{i}\right), \qquad w_{i}\left(0\right)=w_{i,0}, 
\label{follower_dynamic}
\\
&\dot{\tilde{w}}_{k}=\dfrac{1}{N_L}\displaystyle\sum_{h=1}^{N_L}R\left(\tilde{w}_{k},\tilde{w}_{h}\right)\left(\tilde{w}_{h}-\tilde{w}_{k}\right)+u,\qquad\qquad\qquad\qquad\qquad\qquad\qquad \tw_{k}\left(0\right)=\tw_{k,0}, 
\label{leader_dynamic}
\end{align}
where $P(\cdot,\cdot), S(\cdot,\cdot)$ and  $R(\cdot,\cdot)$ are given \emph{compromise functions}, typically taking values in $[0,1]$, measuring the relative importance of the interacting agent in the consensus dynamic. 
The control term $u$ characterizes the strategy of the leaders, and is given by the solution of the following optimal control problem
\begin{equation}\label{eq:minima}
u=\textrm{arg}{\min}\left\{J(u, \underline{w},\underline{\tilde{w}})\right\},\\
\end{equation}
where
\begin{align}\label{J_costfunctional}
&J(u, \underline{w},\underline{\tilde{w}})=\dfrac{1}{2}\int_{0}^{T}\left\{\dfrac{\psi}{N_L}\sum_{h=1}^{N_L}(\tilde{w}_h-w_d)^2+\dfrac{\mu}{N_L}\sum_{h=1}^{N_L}(\tilde{w}_h-m_F)^2\right\}ds+\int_0^{T}{\dfrac{\nu}{2}{u^2}}ds.
\end{align}
In the latter equation $\underline{w}$ and $\underline{\tilde{w}}$ are the vectors with the followers and leaders opinions,
$T$ represents the final time horizon, $w_d$ is the desired opinion and $m_F$ is the average opinion of the followers group at time $t\geq{0}$ defined as
\[
m_F=\dfrac{1}{N_F}\sum_{j=1}^{N_F}w_{j}.
\]
The parameter $\nu>0$, as usual, is a regularization term representing the importance of the control $u$ in the overall dynamic. More precisely, $\nu$ penalizes the action of  the control $u$ in such a way that for large values of $\nu$ the control action vanishes and viceversa.

The problem may also
be formulated as constrained minimization problem for $u^n,\underline{w}^n,\underline{\tw}^n$ in the form
\begin{equation}\label{costrained_min}
\begin{aligned}
&\quad\min J\left(u^n,\underline{w}^n,\underline{\tw}^n\right)\\
&\textrm{subject to }~\eqref{follower_dynamic}-\eqref{leader_dynamic}.
\end{aligned}
\end{equation}
In general the solution of this problems is a difficult task, in particular for nonlinear constrains and non convex functional.
In the following we assume sufficient regularity on the constrains of \eqref{costrained_min}, in such a way that the minimizer fulfills the necessary first order optimality conditions. We refer to \cite{Luen69} for a detailed discussion of necessary and sufficient optimality conditions. 

Thus the control strategy of the leaders' population is based on an interplay of two behaviours weighted by the nonnegative constants $\psi$ and $\mu$ such that $\psi+\mu=1$. On one hand they aim at minimizing the distance with respect to the desired state $w_d$ (\emph{radical behaviour}) and on the other hand they aim at minimizing the distance with respect to the followers' mean opinion (\emph{populistic behaviour}). Therefore, the leaders influence the followers opinion interacting through the function $S(\cdot,\cdot)$ and the followers influence the leaders strategy through their mean opinion in the cost functional (\ref{J_costfunctional}).  

The above optimization problem is approximated using the Boltzmann type optimal control approach recently presented in \cite{AHP13} which corresponds to a binary model predictive control of \eqref{follower_dynamic}-\eqref{eq:minima} in the case of a very large number of agents~\cite{MayneMichalska1990aa, CaBo:04}. 

\subsection*{Instantaneous binary control}
The main idea is to avoid the solution of the dynamic on the whole time interval and to consider a closed-loop strategy for the opinion model in the case of binary interactions. 
Hence, we split the time interval $[0,T]$ in $M$ time intervals of length $\Delta t$ and let  $t^n=\Delta t \; n$ and solve sequentially the optimal control problem in each time interval. This approach is related to the receding horizon strategy, or instantaneous control in the engineering literature, which allows to express the control as a feedback of the state variables. In general, with respect to the associated optimal control problem \eqref{follower_dynamic}-\eqref{eq:minima} this technique furnishes a suboptimal solution. Rigorous results on the properties of $u$ for a constrained quadratic cost functional are discussed, for example, in \cite{MayneMichalska1990aa, CaBo:04}.

More precisely, we approximate both (\ref{follower_dynamic}) and (\ref{leader_dynamic}) by the following discretized binary dynamics 
\begin{equation}\begin{split}\begin{cases}\label{binary_dynamics_follower_new}
w_i^{n+1}=&w^n_i+\alpha P(w^n_i,w^n_j)(w^n_j-w^n_i)+\alpha S(w^n_i,\tw^n_l)(\tw^n_l-w^n_i)\\[+.2cm]
w_j^{n+1}=&w^n_j+\alpha P(w^n_j,w^n_i)(w^n_i-w^n_j)+\alpha S(w^n_j,\tw^n_l)(\tw^n_l-w^n_j)
\end{cases}\end{split}\end{equation}
\begin{equation}\begin{split}\begin{cases}\label{binary_dynamics_leader_new}
\tw_k^{n+1}=&\tw^n_k+\alpha R(\tw^n_k,\tw^n_h)(\tw^n_h-\tw^n_k)+ 2\alpha u^n\\[+.2cm]
\tw_h^{n+1}=&\tw^n_h+\alpha R(\tw^n_h,\tw^n_k)(\tw^n_k-\tw^n_h)+2\alpha u^n
\end{cases}\end{split}\end{equation}
where $\alpha = \Delta t / 2$, $i$ and $j$ are the indexes of the two interacting followers, $l$ the index of an arbitrary leader, $h$ and $k$ the indexes of the two interacting leaders. The control variable $u$ is given by the solution of the following optimization problem
\begin{align}\label{integral_functional}
u^n=\textrm{argmin}~\{ J\left(u^n,\underline{w}^n,\underline{\tw}^n\right)\}
\end{align}
\begin{equation}\label{eq:target}
\begin{aligned}
J\left(u^n,\underline{w}^n,\underline{\tw}^n\right)=\alpha\left(\dfrac{\psi}{2}\sum_{p=\{k,h\}}^{}(\tw^n_{p}-w_d)^2+\dfrac{\mu}{2}\sum_{p=\{k,h\}}^{}(\tw^n_{p}-m_{F}^n)^2+\nu (u^n)^2\right).
\end{aligned}
\end{equation}
In order to solve the minimization problem introduced in \eqref{integral_functional}, we can proceed as in \cite{AHP13} using a standard Lagrange multipliers approach to compute explicitly $u^n$. 
In this way we obtain the feedback control  
\begin{equation}\label{control_s}
2\alpha u^n = -\sum_{p=\{k,h\}}^{}\dfrac{2\alpha^2}{\nu}\left[\psi(\tw_p^{n+1}-w_d)+ \mu(\tw_p^{n+1}-m^{n+1}_{F})\right].
\end{equation}
Note that since the feedback  control $u^n$ in \eqref{control_s} depends on the post interaction opinion the constrained binary interaction \eqref{binary_dynamics_leader_new} is implicitly defined but it can be easily inverted. 
The explicit version of the control reads
\begin{equation}\begin{split}\label{u_partial}
2\alpha u^n=&-\sum_{p=\{k,h\}}\frac{\beta}{2}\left[\psi(\tw^n_p-w_d)+ \mu(\tw^n_p-m^n_{F})\right] -\dfrac{\alpha\beta}{2}(R(\tw^n_k,\tw^n_h)-R(\tw^n_h,\tw^n_k))(\tw^n_h-\tw^n_k),
\end{split}\end{equation}
where we further approximated $m_F^{n+1}$ with $m_F^{n}$ to have a fully explicit expression and introduced the parameter $\beta$ defined as 
\begin{equation}\begin{split}
\beta=&\dfrac{4\alpha^2}{\nu+4\alpha^2}.
\end{split}\end{equation}



\section{Boltzmann type control}
\label{sec:Boltzmann}

In this section, we consider a Boltzmann dynamic corresponding to the above instantaneous control formulation. In order to derive a kinetic equation we introduce a density distribution of followers  $f_F(w,t)$ and leaders $f_L(\tw,t)$ depending on the opinion variables $w,\tw\in I$ and time $t\geq0$.
It is assumed that the followers’ density is normalized to $1$, that
is
 \[
 \int_I f_F(w,t) \, dw = 1,
 \]
whereas  
 \[
 \int_I f_L(\tw,t) \, d\tw = \rho \leq 1.
 \] 
The kinetic model can be derived by considering the change in time of $f_F(w,t)$ and $f_L(\tw,t)$ depending on the interactions with the other individuals and the leaders’ strategy. This change depends on the balance between the gain and loss due to the binary interactions.

\subsection*{Binary constrained interactions dynamic }

Let us consider the pairwise opinions $(w,v)$ and $(\tilde{w},\tilde{v})$, respectively of two followers and two leaders, the corresponding post interaction opinions are computed according with three dynamics, the interaction between two followers, the interaction between follower and leader and finally between two leaders.

The post-interaction opinions $(\tilde{w}^*,\tilde{v}^*)$ of two leaders are given by
\begin{equation}\begin{split}\begin{cases}\label{binary_dynamics_leader0}
\tilde{w}^*=\tilde{w}+\alpha R(\tilde{w},\tilde{v})(\tilde{v}-\tilde{w})+2\alpha u+\tilde{\theta}_{1}\tilde{D}(\tilde{w})\\
\tilde{v}^*=\tilde{v}+\alpha R(\tilde{v},\tilde{w})(\tilde{w}-\tilde{v})+2\alpha u+\tilde{\theta}_{2}\tilde{D}(\tilde{v}),\\
\end{cases}\end{split}\end{equation}  
where the feedback control is defined as
\begin{equation}
\begin{aligned}\label{u_partial}
2\alpha u  =&-\dfrac{\beta}{2}\left[\psi\left((\tilde{w}-w_d)+(\tilde{v}-w_d)\right)+\mu\left((\tilde{w}-m_F)+(\tilde{v}-m_F)\right)\right]\\
&-\dfrac{\alpha\beta}{2}(R(\tilde{w},\tilde{v})-R(\tilde{v},\tilde{w}))(\tilde{v}-\tilde{w}),
\end{aligned}
\end{equation}
and
\begin{equation}
m_F(t)=\int_{I} f_F(w,t)w\,dw.
\end{equation}
Note that the control term is now embedded into the binary interaction and that we considered an additional noise component such that the diffusion variables $\tilde{\theta}_1,\tilde{\theta}_2$ are realizations of a random variable with zero mean and finite variance $\tilde{\sigma}^2$. Moreover the noise influence is weighted by the function $\tilde{D}(\cdot)$, representing the local relevance of diffusion for a given opinion, and such that $0\leq \tilde{D}(\cdot)\leq 1$.

We assume that the opinions $(w^*, v^*)$ in the follower-follower interactions are  derived according to
\begin{equation}\begin{cases}\label{followerfollower_binary}
w^*=w+\alpha P(w,v)(v-w)+\theta_1 D(w), \\
v^*=v+\alpha P(v,w)(w-v)+\theta_2 D(v),
\end{cases}\end{equation}
where the diffusion variables $\theta_1,\theta_2$ are again realizations of a random variable with zero mean, finite variance $\sigma^2$ and $0 \le D(\cdot)\le 1$.
Finally the leader-follower interaction is described for every agents from the leaders’ group, thus in general we have
\begin{equation}\begin{cases}\label{followerleader_binary_1}
w^{**}=w+\alpha S(w,\tilde{v})(\tilde{v}-w)+\hat{\theta}\hat{D}(w) \\
\tilde{v}^{**}=\tilde{v}
\end{cases}\end{equation}

where similar to the previous dynamic, $\hat{\theta}$ is a random variable with zero mean and finite variance $\hat{\sigma}^2$ and $0\le \hat{D}(\cdot)\le 1$



Since we are dealing with a kinetic problem in which the
variable belongs to a bounded domain, namely $I=[-1,1]$, we must deal with 
additional mathematical difficulties in the definition of
agents interactions. In fact, it is essential to consider only
interactions that do not produce values outside the finite
interval. 

For the leaders' interaction if we consider the constrained binary interactions system \eqref{binary_dynamics_leader0}-\eqref{u_partial}, 
without diffusion we obtain that $|\tilde{w}^*-\tilde{v}^*|$ is a contraction if $\alpha \le 1/2$
\begin{eqnarray*}
|\tilde{w}^*-\tilde{v}^*|=
|(\tilde{w}-\tilde{v})-\alpha(\tilde{w}-\tilde{v})(R(\tilde{w},\tilde{v})+R(\tilde{v},\tilde{w}))|
\leq|1-2\alpha| |\tilde{w}-\tilde{v}|.
\end{eqnarray*}
The following proposition gives sufficient conditions to preserve the bounds for the leaders’ interactions (\ref{binary_dynamics_leader0}). 
\begin{prop} \label{prop:v}
Let $r,d_{+}$ and $d_{-}$ be defined as follows
\begin{equation}\label{def:r}
r=\min_{\tilde{v},\tilde{w}\in I}\left[ {R(\tilde{v},\tilde{w})}\right],\qquad
d_{\pm}=\min_{\tilde{w}\in I}\left[{\dfrac{1\mp\tilde{w}}{\tilde{D}(\tilde{w})},\tilde{D}(\tilde{w})\ne 0}\right].
\end{equation}
If $\tilde{v},\tilde{w} \in I$ then $\tilde{v}^*,\tilde{w}^*\in I$ if the following conditions hold
\begin{equation}
\alpha r \ge \dfrac{\beta}{2},
\qquad
d_{-}\left(1-\dfrac{\beta}{2}\right)  \le \tilde{\theta}_i\le d_{+}\left(1-\dfrac{\beta}{2}\right),\quad i=1,2.
\end{equation}
\end{prop}
The proof follows by the same arguments used in \cite{AHP13, DMPW:09} and we omit the details.
On the other hand from the definition of binary interaction between followers  \eqref{followerfollower_binary}, in absence of diffusion, the boundaries are never violated. Indeed since $|w|\le 1$ it follows that $|v-w|\le 1$ and being $0\le P(\cdot,\cdot) \le 1$ it is easily seen that $w^{*},v^{*}\in I$. 

Finally, as shown in \cite{DMPW:09}, the post-interaction opinion of  followers $w^{**}$, in the leader-follower interaction \eqref{followerleader_binary_1}, takes values in the reference interval $I$ if the hypothesis of the following proposition are satisfied. 

\begin{prop}\label{prop:w}
Let $K_{-}$ and $K_+$ be defined as follows
\begin{equation}
K_{\pm}=\min_{{w}\in I}\left[{\dfrac{1\mp{w}}{\hat{D}({w})},\hat{D}({w})\ne 0}\right].
\end{equation}
If $w\in I$ then  $w^{**}\in I$ if the following conditions hold
\begin{equation}
(1-\alpha)K_{-}\le \hat{\theta} \le (1-\alpha)K_+, \quad i=1,2.
\end{equation}
\end{prop}
\subsubsection*{Main properties}
Following the derivation in \cite{PT:13}, for a suitable choice of test functions $\varphi$ we can describe the evolution of $f_F(w,t)$ thanks to the integro-differential equation of Boltzmann type 
\begin{equation}\label{boltzmann_equation}
\dfrac{d}{dt}\int_{I}\varphi(w)f_F(w,t)dw=\left(Q_F(f_F,f_F),\varphi\right)+\left(Q_{FL}(f_L,f_F),\varphi\right)
\end{equation}
where 
\begin{equation}\label{Q_F-psi}
\left(Q_F(f_F,f_F),\varphi\right)=\left<\int_{I^2}B_{int}^{F}(\varphi(w^*)-\varphi(w))f_F(w,t)f_F(v,t)dwdv\right>
\end{equation}
and
\begin{equation}\label{Q_FL-psi}
\left(Q_{FL}(f_F,f_L),\varphi\right)=\left<\int_{I^2}B_{int}^{FL}(\varphi(w^{**})-\varphi(w))f_F(w,t)f_L(\tilde{v},t)dwd\tilde{v}\right>.
\end{equation}
In \eqref{Q_F-psi} and \eqref{Q_FL-psi} we used the notation $\langle\,\cdot\,\rangle$ to indicate the expectations with respect the random variables, respectively $\theta_i, i=1,2$ and $\hat{\theta}$, and the nonnegative interaction kernels $B_{int}^{F},B_{int}^{FL}$ are related to the probability of the microscopic interactions.  
The simplest choice of interaction kernels which guarantees that the post interaction opinions never violate the bounds is given by
\begin{equation}\begin{split}
B^F_{int}=B^F_{int}(w,v,\theta_1,\theta_2)=\eta_F \chi(|w^*|\le 1)\chi(|v^*|\le 1) \\
B^{FL}_{int}=B^{FL}_{int}(w,\tilde{v},\hat{\theta})=\eta_{FL} \chi(|w^{**}|\le 1)\chi(|\hat{v}|\le 1)
\end{split}\end{equation}
where $\eta_F,\eta_{FL}>0$ are constant relaxation rates and $\chi(\cdot)$ is the indicator function.
If we now assume that the interaction parameters are such that $|w^{*}|,|w^{**}|\le 1$ the Boltzmann operators can be written as
\begin{equation}\label{Q_F}
\left(Q_F(f_F,f_F),\varphi\right)=\eta_F\left<\int_{I^2}(\varphi(w^*)-\varphi(w))f_F(w,t)f_F(v,t)dwdv\right>
\end{equation}
\begin{equation}\label{Q_FL}
\left(Q_{FL}(f_F,f_L),\varphi\right)=\eta_{FL}\left< \int_{I^2}(\varphi(w^{**})-\varphi(w))f_F(w,t)f_L(\tilde{v},t)dwd\tilde{v}\right>.
\end{equation}
In order to study the evolution of the average opinion $m_F(t)$, we take $\varphi(w)=w$ in \eqref{boltzmann_equation}. We have that the evolution of the average opinion of followers is
\begin{equation}\begin{split}
\dfrac{d}{dt}m_F(t)=&\dfrac{\eta_F}{2}\left[ \int_{I^2}(w^*+v^*-w-v)f_F(w,t)f_F(v,t)dwdv  \right] \\
&+\eta_{FL}\int_{I^2}(w^{**}-w)f_F(w,t)f_L(\tilde{v},t)dwd\tilde{v},
\end{split}\end{equation}
since the noise in \eqref{followerfollower_binary} has zero mean. From the definition of binary interactions between followers \eqref{followerfollower_binary} and the definition of interaction leader-follower \eqref{followerleader_binary_1} we have 
\begin{equation}\begin{split}\label{m_F}
\dfrac{d}{dt}m_F(t)=&\dfrac{\eta_F}{2}\alpha\int_{I^{2}}(v-w)(P(w,v)-P(v,w))f_F(w,t)f_F(v,t)dwdv\\
&+\eta_{FL}\alpha\int_{I^2}S(w,\tilde{v})(\tilde{v}-w)f_F(w,t)f_L(\tilde{v},t)dwd\tilde{v}.
\end{split}\end{equation}
\begin{rem}\label{remark_PQ}
If we suppose $P$ symmetric, that is $P(w,v)=P(v,w)$, and $S\equiv 1$ we obtain a simplified equation for the time evolution of $m_F$
\begin{equation}\label{mean_follower_evolution_hp}
\dfrac{d}{dt}m_F(t)= \tilde{\eta}_{FL} ~\alpha(m_L(t)-m_F(t))
\end{equation}
where we introduced the notations $\tilde{\eta}_{FL}=\rho~ \eta_{FL}$ and $m_L(t)=\displaystyle \dfrac{1}{\rho}\int_{I} \tilde{w}f_L(\tilde{w},t)d\tilde{w}$.
\end{rem}

The evolution equation for $m_L(t)$ can be found thanks to similar arguments. We can describe the dynamic of $f_L(\tilde{w},t)$ thanks to the following integro-differential equation of Boltzmann type in weak form
\begin{equation}\label{leader_1}
\dfrac{d}{dt}\int_{I}\varphi(\tilde{w})f_L(\tilde{w},t)d\tilde{w}=\left(Q_L(f_L,f_L),\varphi\right)
\end{equation}
where 
\begin{equation}\label{leader_2}
\left(Q_L(f_L,f_L),\varphi\right)=\left<\int_{I^2}B_{int}(\varphi(\tilde{w}^*)-\varphi(\tilde{w}))f_L(\tilde{w},t)f_L(\tilde{v},t)d\tilde{w}d\tilde{v}\right>.
\end{equation}
As before $\langle\,\cdot\,\rangle$ denotes the expectation taken with respect to the random variables $\tilde{\theta_i}, i=1,2$ and $B_{int}$ is related to the probability of the microscopic interactions. A choice which preserves post interaction opinion bounds is 
\begin{equation}
B_{int}=B_{int}(\tilde{w},\tilde{v},\tilde{\theta}_1,\tilde{\theta}_2)=\eta_L\chi(|\tilde{w}^*|\le1)\chi(|\tilde{v}^*|\le1)
\end{equation}
where $\eta_L>0$ is a constant rate and $\chi(\cdot)$ is the indicator function. Let us consider as test function $\varphi(\tilde{w})=\tilde{w}$. Then equation \eqref{leader_1} assumes the form
\begin{equation}
\dfrac{d}{dt}\int_{I}\tilde{w}f_L(\tilde{w},t)d\tilde{w}=\eta_L \left<\int_{I^2}(\tilde{w}^*-\tilde{w})f_L(\tilde{w},t)f_L(\tilde{v},t)d\tilde{w}d\tilde{v}\right>,
\end{equation}
which is equivalent to consider
\[
\dfrac{d}{dt}\int_{I}\tilde{w}f_L(\tilde{w},t)d\tilde{w}=\dfrac{\eta_L}{2}\left< \int_{I^2}(\tilde{v}^*+\tilde{w}^*-\tilde{v}-\tilde{w})f_L(\tilde{w},t)f_L(\tilde{v},t)d\tilde{w}d\tilde{v}\right>.
\]
Then being the noise in \eqref{binary_dynamics_leader0} with zero mean we have
\begin{equation}\begin{split}\label{mean_v_dynamic}
\dfrac{d}{dt}m_L(t)=&\eta_L\alpha\left(1-\beta\right)\dfrac{1}{\rho}\int_{I^2}(R(\tilde{w},\tilde{v})-R(\tilde{v},\tilde{w}))\tilde{v} f_L(\tilde{w},t)f_L(\tilde{v},t) d\tilde{w} d\tilde{v} \\
&+\tilde{\eta}_L\psi\beta(w_d-m_L(t))+\tilde{\eta}_L\beta\mu(m_F(t)-m_L(t)),
\end{split}\end{equation}
where $\tilde{\eta_L}=\rho\eta_L$.
\begin{rem}
If $R(\tilde{w},\tilde{v})=R(\tilde{v},\tilde{w})$ equation \eqref{mean_v_dynamic} becomes
\begin{equation}
\dfrac{d}{dt}m_L(t)=\tilde{\eta}_L\psi\beta(w_d-m_L(t))+\tilde{\eta}_L\mu\beta(m_F(t)-m_L(t)).
\end{equation}
Moreover if the assumptions on $P$ and $S$ in Remark \ref{remark_PQ} hold we obtain the following closed system of differential equations for the mean opinions $m_L$ and $m_F$
\begin{equation}\begin{cases}\label{mean_system_hp}
\displaystyle\dfrac{d}{dt}m_L(t)=\tilde{\eta}_L\psi\beta(w_d-m_L(t))+\tilde{\eta}_L \mu \beta(m_F(t)-m_L(t))\\[+.25cm]
\displaystyle\dfrac{d}{dt}m_F(t)=\tilde{\eta}_{FL} \alpha (m_L(t)-m_F(t)).
\end{cases}\end{equation}
Straightforward computations show that the exact solution of the above system has the following structure
\begin{equation}\begin{cases}\label{mean_soultion_hp}
m_L(t)=  C_1\exp\left\{-|\lambda_1| t\right\}+C_2\exp\left\{-|\lambda_2| t\right\}+w_d\\[+.25cm]
m_F(t)= \displaystyle C_1\left(1+\frac{\lambda_1}{\beta\mu\tilde{\eta}_L}\right)\exp\left\{-|\lambda_1| t\right\}+C_2\left(1+\frac{\lambda_2}{\beta\mu\tilde{\eta}_L}\right)\exp\left\{-|\lambda_2| t\right\}+w_d\\
\end{cases}\end{equation}
where $C_1,C_2$ depend on the initial data $m_F(0),m_L(0)$ in the following way
\begin{align}
C_1=&-\frac{1}{\lambda_1-\lambda_2}\left( (\beta\tilde{\eta}_Lm_L(0)  + \lambda_2)m_L(0) - \mu\beta\tilde{\eta}_L m_F(0)-(\lambda_2+\beta\tilde{\eta}_L\psi)w_d\right)\\
C_2=& \quad\frac{1}{\lambda_1-\lambda_2}\left((\beta\tilde{\eta}_Lm_L(0)  + \lambda_1)m_L(0) - \mu\beta\tilde{\eta}_L m_F(0)-(\lambda_1+\beta\tilde{\eta}_L\psi)w_d\right)\end{align}
with 
\begin{equation}
\lambda_{1,2}=-\frac{1}{2}\left(\alpha\tilde{\eta}_{FL}+\beta\tilde{\eta}_{L}\right)\pm\frac{1}{2}\sqrt{(\alpha\tilde{\eta}_{FL}+\beta\tilde{\eta}_{L})^2-4\psi\alpha\beta\tilde{\eta}_L\tilde{\eta}_{FL}}.
\end{equation}
Note that $\lambda_{1,2}$ are always negative, this assures that the contribution of the initial averages, $m_L(0),m_F(0)$, vanishes as soon as time increases and the mean opinions of leaders and followers converge towards the desired state $w_d$.
\end{rem}
We now take into account the evolution of the second order moments
\[
E_F(t)=\int_I w^2 f_F(w,t)~dw,\qquad E_L(t)=\dfrac{1}{\rho}\int_I \tilde{w}^2 f_L(\tilde{w},t)~d\tilde{w}.
\]
First we analyze the followers group from equation \eqref{boltzmann_equation} with test functions $\varphi(w)=w^2$, we have
\begin{equation}\begin{aligned}\label{variance_general}
\dfrac{d}{dt}E_F(t)=&\dfrac{\eta_F}{2}\left<\int_{I^2}\left((w^*)^2+(v^*)^2-w^2-v^2\right)f_F(w,t)f_F(v,t)dwdv\right>\\
&+\eta_{FL}\left<\int_{I^2}\left((w^{**})^2-w^2\right)f_F(w,t)f_L(\tilde{v},t)dwd\tilde{v}\right>.
\end{aligned}\end{equation}
Thanks to \eqref{followerfollower_binary} -\eqref{followerleader_binary_1}, in the simplified case $P\equiv S \equiv 1$, we obtain
\begin{equation}\begin{split}\label{energy_followers}
\dfrac{d}{dt}E_F(t)=&2\eta_{F}\alpha(\alpha-1)(E_F(t)-m_F^2(t))+\tilde{\eta}_{FL} \alpha^2(E_L+E_F-2m_L(t)m_F(t))\\
&+2\alpha\tilde{\eta}_{FL}(m_F(t)m_L(t)-E_F(t))+\eta_{F}\sigma^2\int_ID^2(w)f_F(w,t)dw\\         &+\tilde{\eta}_{FL}\hat{\sigma}^2 \int_I\hat{D}^2(w)f_F(w,t)dw.
\end{split}\end{equation}
Finally, for the leaders group let us consider the function $\varphi(\tilde{w})=\tilde{w}^2$ in \eqref{leader_1} and the case $R\equiv1$. Then thanks to equation \eqref{binary_dynamics_leader0}  we obtain
\begin{equation}
\begin{aligned}\label{energy_leaders}
\dfrac{d}{dt}E_L(t)=& \dfrac{\eta_L}{2}\dfrac{1}{\rho}\left<\int_{I^2}((\tilde{w}^*)^2+(\tilde{v}^*)^2-\tilde{w}^2-\tilde{v}^2)f_L(\tilde{w},t)f_L(\tilde{v},t)d\tilde{w}d\tilde{v}\right> \\
= & \tilde{\eta}_L \Big[2\alpha(\alpha-1)(E_L(t)-m_L^2(t))-\dfrac{\beta}{2}(2-\beta)(E_L(t)+m_L^2(t))\\
&+2\beta(1-\beta)(\psi w_d+\mu m_F(t))m_L(t)+\beta^2(\psi w_d\\
&+\mu m_F(t))^2+\tilde{\sigma}^2\int_I\tilde{D}^2(\tilde{w})f_L(\tilde{w},t)d\tilde{w} \Big].
\end{aligned}
\end{equation}
In absence of diffusion, since $m_F(t),m_L(t)\rightarrow w_d$ as $t\rightarrow \infty$, it follows that $E_F(t),E_L(t)$ converge toward $w_d^2$. Then the quantities 
\begin{equation}\begin{split}
\int_I f_F(w,t) (w-w_d)^2dw=&E_F(t)+w_d^2-2m_F(t)w_d \\
\dfrac{1}{\rho}\int_I f_L(\tilde{w},t)(\tilde{w}-w_d)^2d\tilde{w}=&E_L(t)+w_d^2-2m_L(t)w_d 
\end{split}
\end{equation}
go to zero as $t\rightarrow \infty$, i.e. under the above assumptions the steady state solutions have the form of a Dirac delta centered in the target opinion $w_d$.

\section{Fokker-Planck Modeling}\label{sec:FP}
In the general case,  it is quite difficult to obtain analytic results on the large time behaviour of the kinetic equation \eqref{boltzmann_equation}. 
A step towards the simplification of the analysis, is the derivation of asymptotic states of the Boltzmann model resulting in simplified Fokker-Planck type models, for which the study of the asymptotic properties is easier \cite{PT:13}. 
In order to obtain such simplification we will follow the approach usually referred as {quasi--invariant opinion} limit~\cite{PT:13,Tos06}, which is closely related to the so-called {grazing collision limit} of the Boltzmann equation (see \cite{FPTT12, Vill98}). 

\subsection*{Quasi invariant opinion limit} 
The main idea is to rescale the interaction frequencies  $\eta_L,\eta_F,\eta_{FL}$, the propensity strength $\alpha$, the diffusion variances $\tilde{\sigma}^2,\sigma^2,\hat{\sigma}^2$ and the action of the control $\nu$ at the same time, in order to maintain, at level of the asymptotic procedure, memory of the microscopic interactions.

Let us introduce the parameter $\varepsilon>0$, and consider the rescaling 
\begin{equation}\label{eq:scaling}
\begin{aligned}
&\alpha=\varepsilon,\qquad \nu=\varepsilon \kappa, \qquad \sigma^2=\varepsilon \varsigma^2, \qquad \hat{\sigma}^2=\varepsilon\hat{\varsigma}^2,\qquad \tilde{\sigma}^2=\varepsilon \tilde{\varsigma}^2,\\
&\eta_F=\dfrac{1}{c_F\varepsilon},\qquad \eta_{FL}=\dfrac{1}{c_{FL}\varepsilon},\qquad \eta_L=\dfrac{1}{c_L\varepsilon}, \qquad  \beta=\dfrac{4\varepsilon}{\kappa+4\varepsilon}.
\end{aligned}
\end{equation}
This corresponds to the situation where the interaction operator concentrates on binary interactions which produce a very small change in the opinion of the agents. From a modeling viewpoint, we require that the scaling \eqref{eq:scaling} in the limit $\varepsilon\to 0$ preserves the main macroscopic properties of the kinetic system. To this extent, let us consider the evolution of the scaled first two moments under the simplifying hypothesis $P,R$ symmetric  and $S\equiv 1$.

The evolution of the mean opinions described in the system \eqref{mean_system_hp} rescales as
\begin{equation}\begin{cases}\begin{split}
\dfrac{d}{dt}m_F(t)=&\varepsilon \dfrac{1}{c_{FL}\varepsilon}(m_L(t)-m_F(t)) \\
\dfrac{d}{dt}m_L(t)=&\dfrac{\psi}{c_L\varepsilon}\dfrac{4\varepsilon}{\kappa+4\varepsilon}(w_d-m_L(t))+\dfrac{\mu}{c_L\varepsilon}\dfrac{4\varepsilon}{\kappa+4\varepsilon}(m_F(t)-m_L(t))
\end{split}\end{cases}\end{equation}
which as $\varepsilon \rightarrow 0$ yields
\begin{equation}\begin{cases}\label{system_mean}
\vspace{5 mm} \dfrac{d}{dt}m_F(t)=\dfrac{\rho}{c_{FL}}(m_L(t)-m_F(t)) \\ 
\dfrac{d}{dt}m_L(t)=\dfrac{4\rho}{c_L\kappa}\left[\psi(w_d-m_L(t))+\mu(m_F(t)-m_L(t)))\right].
\end{cases}\end{equation}
The second moment equations \eqref{energy_followers} and \eqref{energy_leaders} are then scaled as follows
\begin{equation}\begin{split}\label{energy_followers_scaled}
\dfrac{d}{dt}E_F(t)=&(\varepsilon-1)\dfrac{2}{c_F}(E_F(t)-m_F^2(t))+\dfrac{\varepsilon \rho}{c_{FL}}(E_L(t)+E_F(t)-2m_L(t)m_F(t))\\
&+\dfrac{2\rho}{c_{FL}}(m_F(t)m_L(t)-E_F(t))+\dfrac{\varsigma^2}{c_F}\int_I{D^2(w)}f_F(w,t)dw\\
&+\dfrac{\hat{\varsigma}^2\rho}{c_{FL}}\int_I{\hat{D}^2(w)}f_F(w,t)dw
\end{split}\end{equation}
\begin{equation}\begin{split}
\dfrac{d}{dt}E_L(t)=&\dfrac{\rho}{c_L\varepsilon}\left[ 2\varepsilon(\varepsilon-1)(E_L(t)-m_L^2(t))-\dfrac{2\varepsilon}{\kappa+4\varepsilon}(2-\dfrac{4\varepsilon}{\kappa+4\varepsilon})(E_L(t)+m_L^2(t))\right. \\
&+\left.\dfrac{8\varepsilon}{\kappa+4\varepsilon}(1-\dfrac{4\varepsilon}{\kappa+4\varepsilon})(\psi w_d +\mu m_F(t))m_L(t)\right. \\
&\left. +\left(\dfrac{4\varepsilon}{\kappa+4\varepsilon}\right)^2(\psi w_d + \mu m_F(t))^2+\tilde{\sigma}^2\int_I \tilde{D}^2(\tilde{w})f_L(\tilde{w},t)d\tilde{w}.\right]
\end{split}\end{equation}
and as $\varepsilon \rightarrow 0$ 
%
%
we obtain  
\begin{equation}
\begin{aligned}
\dfrac{d}{dt}E_F(t)=&-\dfrac{2}{c_F}(E_F(t)-m_F^2(t))+\dfrac{2\rho}{c_{FL}}(m_F(t)m_L(t)-E_F(t))\\    
&+\dfrac{\varsigma^2}{c_F}\int_I{D^2(w)}f_F(w,t)dw+\dfrac{\hat{\varsigma}^2\rho}{c_{FL}}\int_I{\hat{D}^2(w)}f_F(w,t)dw \\
\dfrac{d}{dt}E_L(t)=&-\dfrac{2\rho}{c_L}(E_L(t)-m_L^2(t))-\dfrac{4\rho}{c_L \kappa}(E_L(t)+m_L^2(t))\\
&+ \dfrac{8\rho}{c_L \kappa }(\psi w_d+\mu m_F(t))m_L(t)+\dfrac{\tilde{\varsigma}^2\rho}{c_L}\int_I \tilde{D}^2(\tilde{w})f_L(\tilde{w},t)d\tilde{w}.
\end{aligned}
\end{equation}
Therefore the asymptotic scaling preserve the behaviour of the first two moments of the solution. We show how this approach leads to a constrained Fokker--Planck system for the description of the opinion distribution of leaders and followers. We present formal computation, following the same arguments in \cite{PT:13,Tos06} it is possible to give a rigorous mathematical basis of our derivation. Here we omit the details for the sake of brevity.


\subsection*{Fokker-Plank equations}
The scaled equation \eqref{boltzmann_equation} reads
\begin{equation}\begin{split}\label{boltzmann_scaled}
\dfrac{d}{dt}\int_{I}\varphi(w)f_F(w,t)dw=&\dfrac{1}{c_F\varepsilon}\left<\int_{I^2}(\varphi(w^*)-\varphi(w))f_F(w,t)f_F(v,t)dwdv\right>\\
&+\dfrac{1}{c_{FL}\varepsilon}\left<\int_{I^2}(\varphi(w^{**})-\varphi(w))f_F(w,t)f_L(\tilde{w},t)dwd\tilde{w}\right>.
\end{split}\end{equation}
Considering the second order Taylor expansion of $\varphi$ around $w$ we obtain
\begin{equation}\begin{split}\label{taylor_exp}
\varphi(w^*)-\varphi(w)=(w^*-w)\varphi '(w)+\dfrac{1}{2}(w^*-w)^2\varphi ''(\bar{w}) \\
\varphi(w^{**})-\varphi(w)=(w^{**}-w)\varphi '(w)+\dfrac{1}{2}(w^{**}-w)^2\varphi ''(\hat{w}) 
\end{split}\end{equation}
where for some $0 \le \vartheta_1,\vartheta_2 \le 1$
\begin{equation*}
\bar{w}=\vartheta_1 w^*+(1-\vartheta_1)w,\qquad \hat{w}=\vartheta_2 w^{**}+(1-\vartheta_2)w.
\end{equation*}
%
Taking into account the binary interactions \eqref{followerfollower_binary}-\eqref{followerleader_binary_1} in \eqref{taylor_exp}, and  substituting in \eqref{boltzmann_scaled}, we obtain a second order approximation of the dynamic. In the limit $\varepsilon \rightarrow 0$ the leading order is given by
\begin{equation}\begin{aligned}
\dfrac{d}{dt}\int_I \varphi(w)f_F(w)dw=&\dfrac{1}{c_F}\left[\int_{I^2}P(w,v)(v-w)\varphi '(w)f_F(w,t)f_F(v)dwdv \right] \\
& +\dfrac{1}{c_{FL}}\left[\int_{I^2}S(w,\tilde{w})(\tilde{w}-w)\varphi '(w)f_F(w)f_L(\tilde{w})dwd\tilde{w} \right] \\
& +\dfrac{1}{2}\dfrac{\varsigma^2}{c_F}\int_I \varphi''(w)D^2(w)f_F(w,t)dw\\
&+\dfrac{1}{2}\dfrac{\hat{\varsigma}^2\rho}{c_{FL}}\int_I \varphi''(w) \hat{D}^2(w)f_F(w,t)dw.
 \end{aligned}\end{equation} 
Integrating back by parts the last expression we obtain the Fokker-Planck equation for the followers' opinion distribution
\begin{equation}\label{fokker-planck-follower}
\begin{aligned}
\dfrac{\partial f_F}{\partial{t}}&+\dfrac{\partial}{\partial w}\left(\dfrac{1}{c_F}K_{F}[f_F](w)+\dfrac{1}{c_{FL}}K_{FL}[f_L](w)\right)f_F(w)=&
\dfrac{1}{2}\dfrac{\partial ^2}{\partial \tilde{w}^2}\left(\dfrac{\varsigma^2}{c_F}\tilde{D}^2(\tilde{w})+\dfrac{\hat{\varsigma}^2\rho}{c_{FL}}\right)f_F({w}),
\end{aligned}
\end{equation}
where 
\begin{equation}
K_{F}[f_F](w)=\int_I P(w,v)(v-w)f_F(v,t)dv,
\qquad K_{FL}[f_L](w)=\int_I S(w,\tilde{w})(\tilde{w}-w)f_L(\tilde{w})d\tilde{w}.
\end{equation}
Following the same strategy we obtain the analogous result for the leaders' opinion distribution 
\begin{equation}\label{fokker-planck-leader}
\begin{aligned}
\dfrac{\partial f_L}{\partial t} + \dfrac{\partial}{\partial \tilde{w}}\left(\dfrac{\rho}{c_L}H[f_L](\tilde{w})+\dfrac{1}{c_L}K_L[f_L](\tilde{w})\right)f_L(\tilde{w})=\dfrac{1}{2}\dfrac{\tilde{\varsigma}^2\rho}{c_L}\dfrac{\partial ^2}{\partial \tilde{w}^2}\tilde{D}^2(\tilde{w})f_L(\tilde{w})
\end{aligned}
\end{equation}
where
\begin{equation} 
K[f_L](\tilde{w})=\int_I R(\tilde{w},\tilde{v})(\tilde{v}-\tilde{w})f_L(\tilde{v},t)d\tilde{v} 
\end{equation}
and
\begin{equation}\begin{split}
H[f_L](\tilde{w})=&\dfrac{2\psi}{\kappa}\left(\tilde{w}+m_L(t)-2w_d\right)+\dfrac{2\mu}{\kappa}\left(\tilde{w}+m_L(t)-2m_F(t)\right).
\end{split}\end{equation}

\subsection*{Steady state solutions}
In this section we show that in some cases it is possible to find explicit stationary states of the Fokker-Planck system of equation described in \eqref{fokker-planck-follower} and \eqref{fokker-planck-leader}. Here we restrict to the simplified situation where every interaction function is constant and unitary, i.e. $P\equiv S\equiv R\equiv 1$, and 
\begin{equation}\begin{split}
D(w)=\tilde{D}(w)=\hat{D}(w)=1-w^2.
\label{eq:cdiff}
\end{split}\end{equation}
The steady state of equations \eqref{fokker-planck-follower} and \eqref{fokker-planck-leader} is solution of the following equations
\begin{equation}
\begin{aligned}
&\left(\dfrac{1}{c_F}(m_F-w)+\dfrac{\rho}{c_{FL}}(m_L-w)\right)f_{F,\infty}=\dfrac{1}{2}\left(\dfrac{\varsigma^2}{c_F}+\dfrac{\hat{\varsigma}^2\rho}{c_{FL}}\right)\dfrac{\partial}{\partial w}D^2(w)f_{F,\infty},\\
&\left(\dfrac{2\psi}{\kappa}[\tilde{w}-2w_d-m_L]+\dfrac{2\mu}{\kappa}[\tilde{w}-2m_F+m_L]\right)f_{L,\infty}=\dfrac{1}{2}\dfrac{\tilde{\varsigma}^2\rho}{c_L}\dfrac{\partial}{\partial \tilde{w}}D^2(\tilde{w})f_{L,\infty}.
\end{aligned}
\end{equation}
As soon as $t\rightarrow \infty$, thanks to equation \eqref{system_mean}, the followers and the leaders’ mean opinion $m_F$ and $m_L$ relax to the desired opinion $w_d$. Then  
\begin{equation}
\left(\dfrac{1}{c_F}+\dfrac{\rho}{c_{FL}}\right)(w_d-w)f_{F,\infty}=\dfrac{1}{2}\left(\dfrac{\varsigma^2}{c_F}+\dfrac{\hat{\varsigma}^2\rho}{c_{FL}}\right)\dfrac{\partial}{\partial w}D^2(w)f_{F,\infty},
\end{equation}
that is 
\begin{equation}
\left(\dfrac{1}{c_F}+\dfrac{\rho}{c_{FL}}\right)(w_d-w)\dfrac{g_{F}}{D^2(w)}=\dfrac{1}{2}\left(\dfrac{\varsigma^2}{c_F}+\dfrac{\hat{\varsigma}^2\rho}{c_{FL}}\right)\dfrac{\partial}{\partial w}g_{F}
\end{equation}
where $g_F=D^2(w)f_{F,\infty}$. This implies
\begin{equation}
g_{F,\infty}=a_F\exp\left\{{-\dfrac{2}{b_F}\int_0^{w}\dfrac{z-w_d}{(1-z^2)^2}}dz\right\},\qquad 
b_F=\dfrac{\varsigma^2 c_{FL}+\hat{\varsigma}^2 c_F\rho}{c_{FL}+c_F\rho},
\end{equation}
and $a_F$ is a normalization constant such that $\int_I g_{F,\infty}dw=1$. Finally we have
\begin{equation}
f_{F,\infty}=\dfrac{a_F}{(1-w^2)^2}\exp\left\{{-\dfrac{2}{b_F}\int_0^{w}\dfrac{z-w_d}{(1-z^2)^2}}dz\right\}.
\end{equation}
Similarly we can find the steady state $f_{L,\infty}$ as a solution of the equation 
\begin{equation}\label{steady_leader}
-\left(\dfrac{2\psi}{\kappa}+\dfrac{2\mu}{\kappa}\right)(w_d-\tilde{w})\dfrac{g_{L,\infty}}{D^2(\tilde{w})}=\dfrac{1}{2}\dfrac{\tilde{\varsigma}^2\rho}{c_L}\dfrac{\partial}{\partial \tilde{w}}g_{L,\infty},
\end{equation}
where $g_{L,\infty}=f_{L,\infty}D^2(w)$. The solution of the differential equation \eqref{steady_leader} is given by
\begin{equation}
g_{L,\infty}=a_L \exp\left\{-\dfrac{2}{b_L}\int_0^{\tilde{w}}\left(\dfrac{z-w_d}{(1-z^2)^2}\right)dz\right\},\qquad
b_L=\dfrac{\bar\varsigma\rho \kappa}{2 c_L (\psi+\mu)},
\end{equation}
and $a_L$ is chosen such that the mass of $g_{L,\infty}$ is equal to $\rho$. Then the steady state is 
\begin{equation}
f_{L,\infty}=\dfrac{a_L}{(1-\tilde{w}^2)^2}\exp\left\{-\dfrac{2}{b_L}\int_0^{\tilde{w}}\left(\dfrac{z-w_d}{(1-z^2)^2}\right)dz\right\}.
\end{equation}

\section{Numerical Simulation}\label{sec:Numeric}
In this section we present several numerical results concerning the numerical simulation of the Boltzmann type control model introduced in the previous paragraphs. All the results have been computed by a Monte Carlo method for the Boltzmann model (see \cite{PT:13} for more details) in the Fokker-Planck regime $\varepsilon=0.01$ under the scaling (\ref{eq:scaling}). 
In the numerical tests we assume that the five per cent of the population is composed by opinion leaders, see for example \cite{DMPW:09}. Note that, for clarity, in all figures the leaders’ profiles have been magnified by a factor $10$. The regularization term in the control is fixed to $\nu=1$. The random diffusion effects have been computed in the case (\ref{eq:cdiff}) for a uniform random variable with scaled variance $\varsigma^2=\tilde\varsigma^2=\hat\varsigma^2=0.01$. It is easy to check that the above choices preserve the bounds in the numerical simulations. 
First we present some test cases with a single population of leaders as discussed in our theoretical analysis. Then we consider the case of multiple leaders’ populations with different time-dependent strategies. This leads to more realistic applications of our arguments, introducing the concept of competition between leaders’ populations. For the sake of simplicity we fix constant interaction functions $P(\cdot,\cdot)\equiv 1$ and $R(\cdot,\cdot)\equiv 1$ and the remaining scaled computational parameters have been summarized in Table \ref{tab:par}. 
\begin{table}[h]
\centering
\caption{Computational parameters for the different test cases.}
\begin{tabular}{c|c|c|ccccc|ccccc}
\hline
\hline
Test & $S(\cdot,\cdot)$ & $c_F$& $\hat{c}_{FL}$&$\hat{c}_{L}$ & $\rho$ &$\psi$&$w_d$ \\
\hline
\hline
1a& 1& 1& 0.1&0.1 & 0.05& 0.5& 0.5 & & & & &\\
\hline
1b& eq. (\ref{eq:bc}) & 1& 0.1&0.1 & 0.05& 0.5& 0.5 & & & & &\\
\hline
\hline
~& $S(\cdot,\cdot)$ & $c_F$& $\hat{c}_{FL_1}$&$\hat{c}_{L_1}$ & $\rho_1$ &$\psi_1$&$w_{d_1}$& $\hat{c}_{FL_2}$&$\hat{c}_{L_2}$ & $\rho_2$ &$\psi_2$&$w_{d_2}$\\
\hline
\hline
2& 1 & 1& 0.1&0.1 & 0.05& 0.5& 0.5 & 0.1&0.1 & 0.05& 0.5& -0.5\\
\hline
3& 1 & 1& 0.1&0.1 & 0.05& eq. (\ref{psi_timedep}) & 0.5 & 1&0.1 & 0.05& eq. (\ref{psi_timedep})&-0.5\\
\hline
\end{tabular}
\label{tab:par}
\end{table}

\subsection*{Test 1. Leaders driving followers}
\begin{figure}[t]
\centering 
\includegraphics[scale=.38]{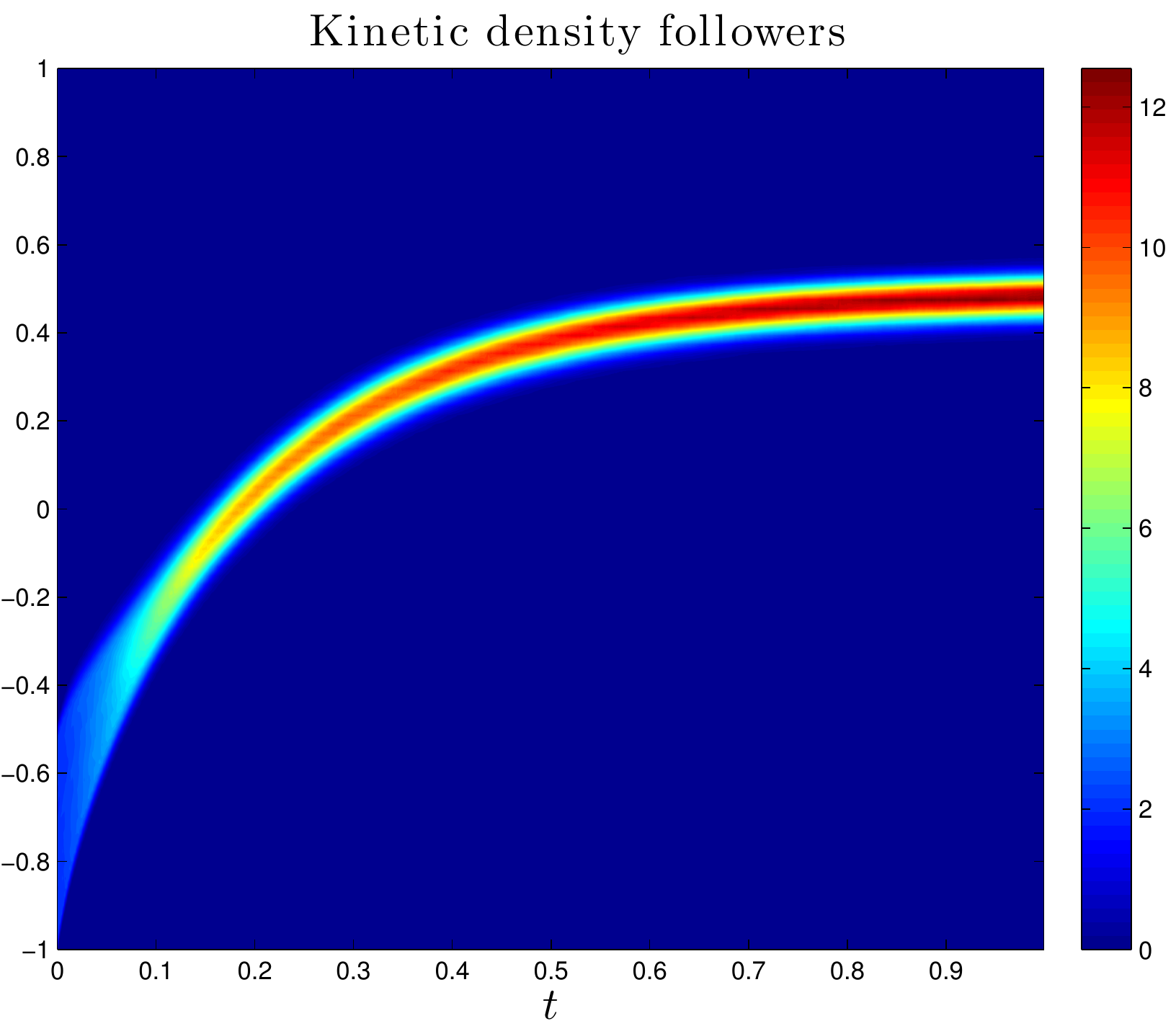}
\includegraphics[scale=.38]{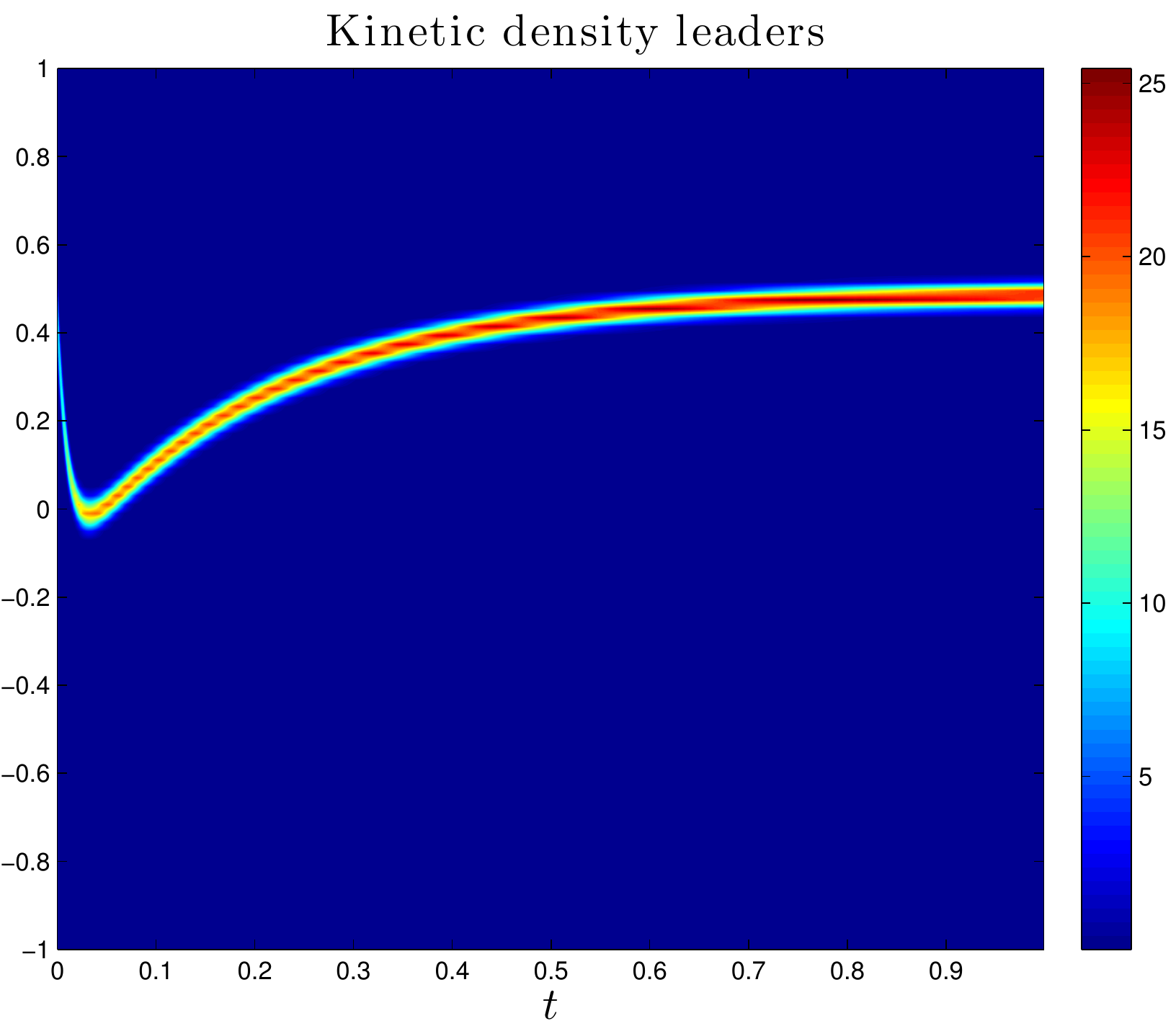}
\caption{Test 1a: Kinetic densities evolution over the time interval $[0,1]$ for a single  population of leaders.}
\label{fig:kinetic_dens_2pop}
\end{figure}

In the first test case we consider a single population of leaders driving followers described by the following system of Boltzmann equations
\begin{equation}\begin{cases} \vspace{2 mm}
\dfrac{d}{dt}\displaystyle\int_I\varphi(w)f_F(w,t)dw=\left(Q_F(f_F,f_F),\varphi \right)+\left(Q_{FL}(f_F,f_L),\varphi\right) &\\
\dfrac{d}{dt}\displaystyle\int_I\varphi(\tilde{w})f_L(\tilde{w},t)d\tilde{w}=\left(Q_L(f_L,f_L),\varphi \right).&
\end{cases}\end{equation}
Numerical experiments show that the optimal control problem is capable to introduce a non monotone behaviour of $m_L(t)$. 
We report the evolution, over the time interval $[0, 1]$, of the the kinetic densities $f_F(w,t)$ and $f_L(\tilde{w},t)$ in Figure \ref{fig:kinetic_dens_2pop} for constant interaction functions $P,R$ and $S$. The initial distributions $f_F\sim U([-1,-0.5])$ and $f_L\sim N(w_d,0.05)$ where $U(\cdot)$ and $N(\cdot,\cdot)$ denote, as usual, the uniform and the normal distributions. We used the compact notations 
\be
\hat{c}_{FL}={c_{FL}}/{\rho},\qquad \hat{c}_L={c_L}/{\rho}.
\ee
This non monotone behaviour shows that the leaders use a combination of populistic and radical strategy to drive the followers towards their desired state. In an electoral context, this is a characteristic which can be found in populistic radical parties, which typically include non-populist ideas and their leadership generates through a dense network of radical movements \cite{Cas07}. 

\begin{figure}[t]
\centering
\includegraphics[scale=.38]{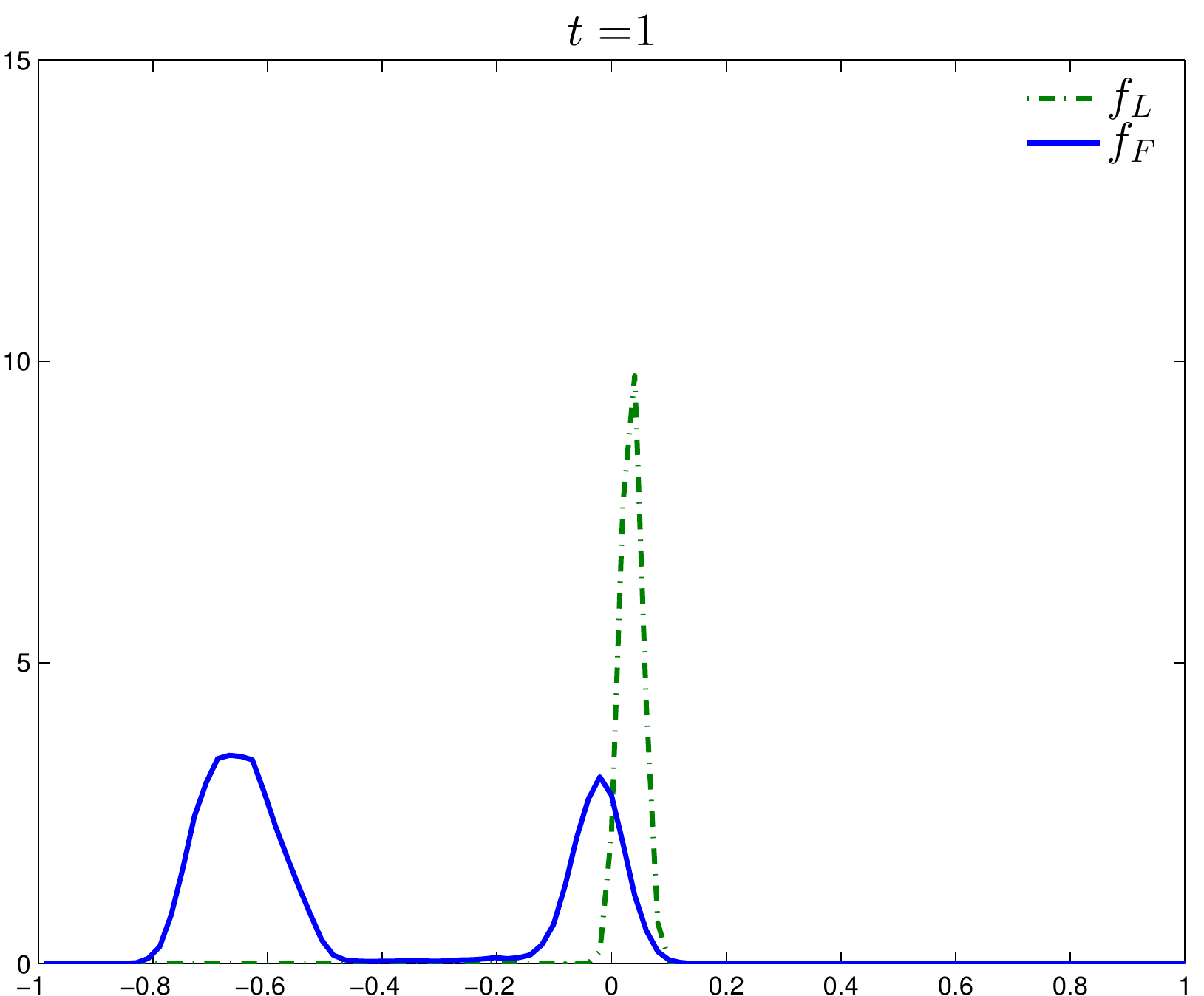}
\includegraphics[scale=.38]{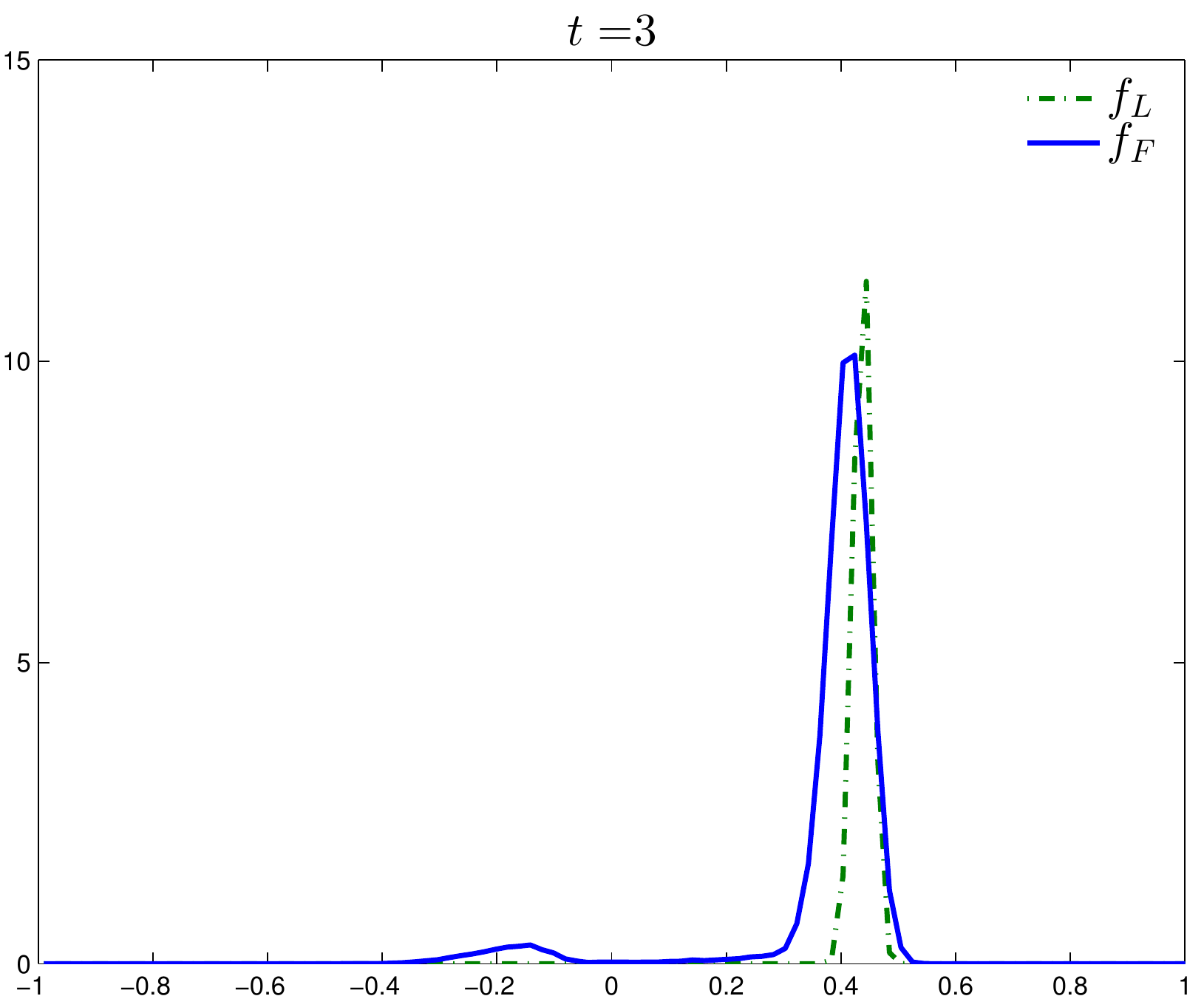}
\caption{Test 1b: Kinetic densities at different times for a single  population of leaders with bounded confidence interaction.}
\label{fg:dc1}
\end{figure}

\begin{figure}[t]
\centering 
\includegraphics[scale=.38]{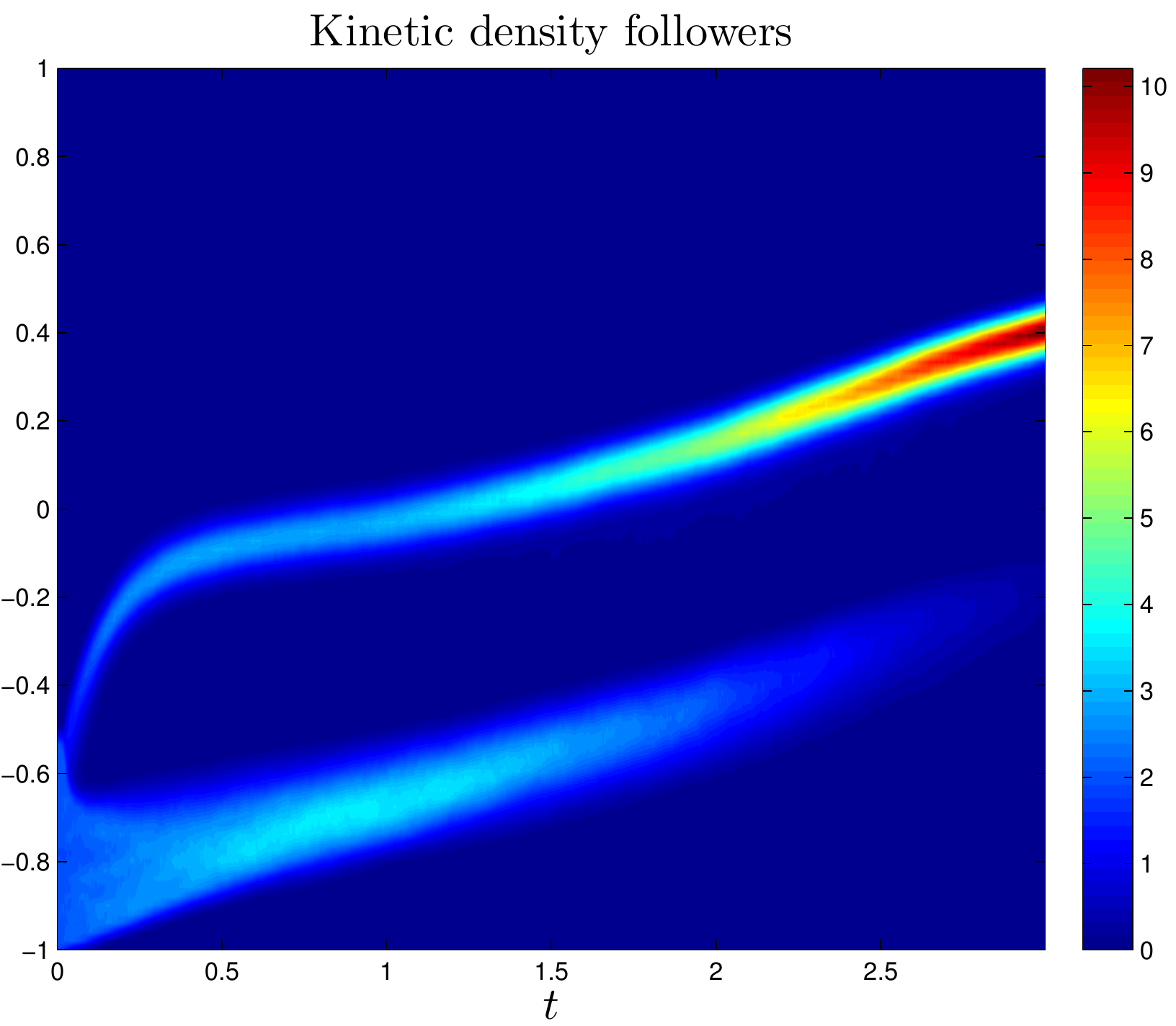}
\includegraphics[scale=.38]{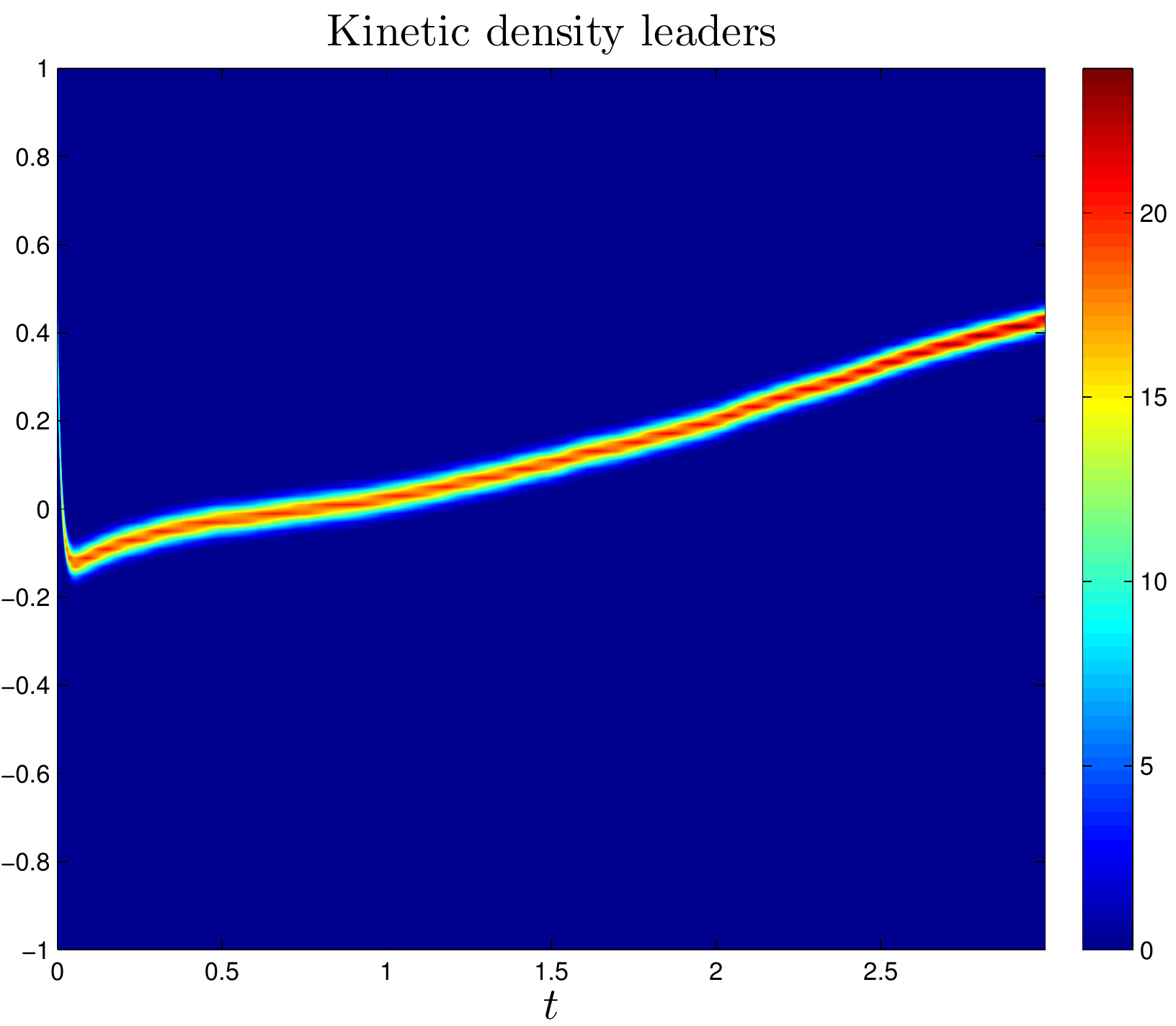}
\caption{Test 1b: Kinetic densities evolution over the time interval $[0,3]$ for a single  population of leaders with bounded confidence interaction.}
\label{fg:dc2}
\end{figure}

Next we consider a bounded confidence model for the leader-follower interaction with
\be
S(w,\tilde{w})=\chi(|w-\tilde{w}|\leq \Delta),
\label{eq:bc}
\ee
where $0 \leq \Delta\leq 2$. In the simulation we assume $\Delta = 0.5$ and use the same initial data of the previous case. It is interesting to observe how the model is capable to reproduce a realistic behavior where the leaders first are able to attract a small group of followers which subsequently are capable to drive the whole majority towards the desired state (see Figures \ref{fg:dc1} and \ref{fg:dc2}). 

\begin{figure}[p]
\centering
\includegraphics[scale=.38]{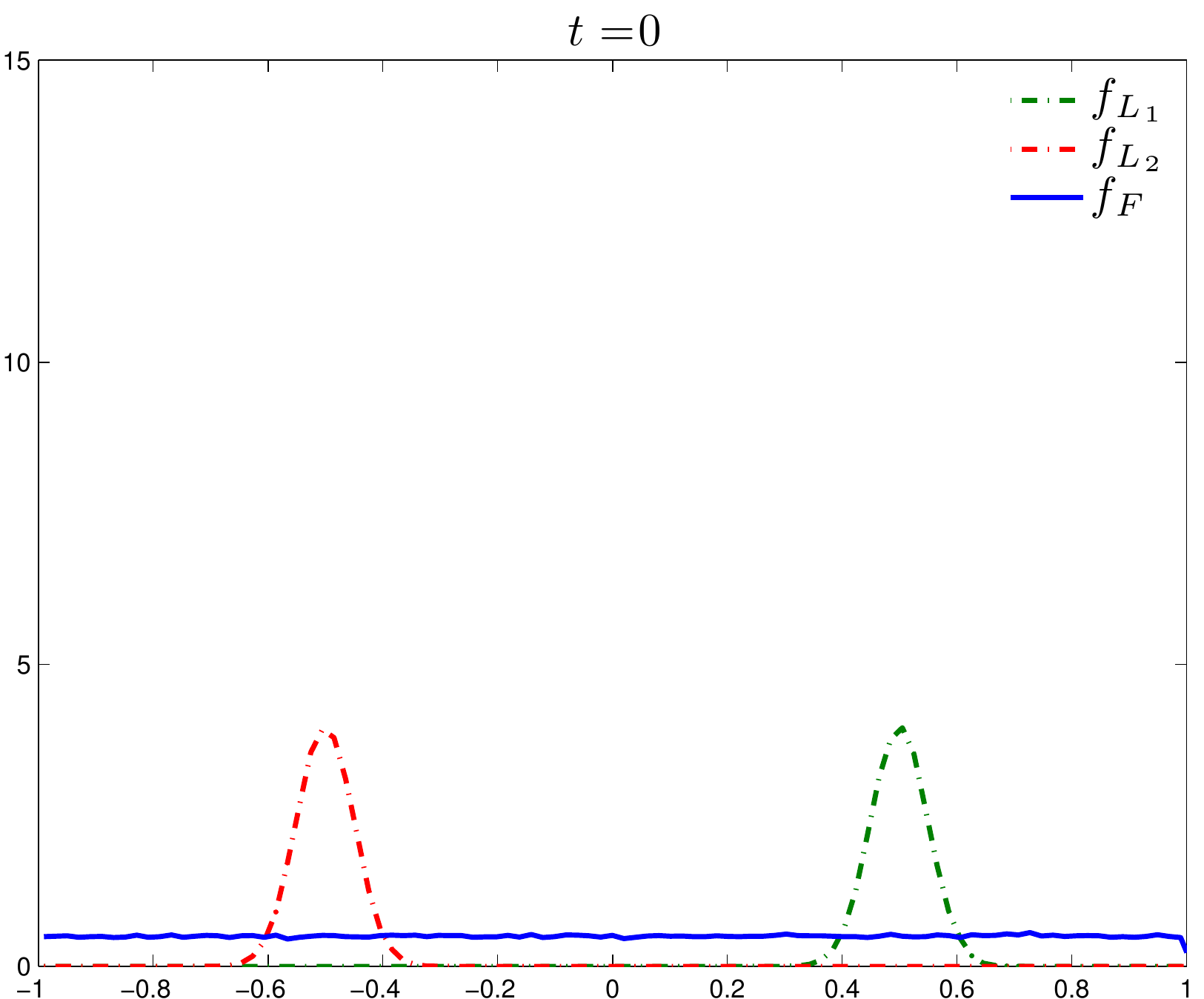}
\includegraphics[scale=.38]{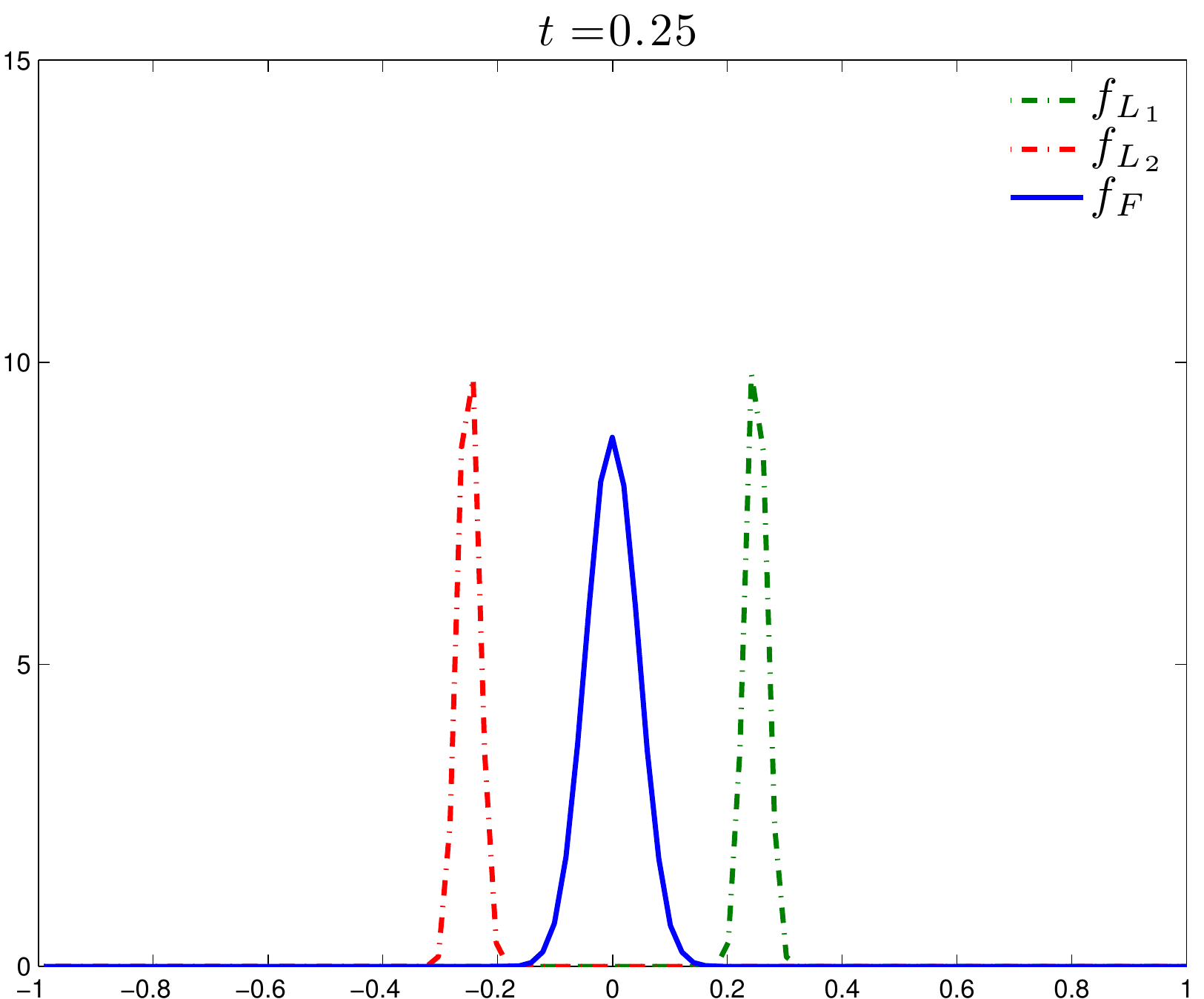}
\caption{Test 2: Kinetic densities at different times reproducing a Hotelling-like model behaviour for two populations of leaders.}
\label{fig:hot1}
\end{figure}

\begin{figure}[p]
\centering
\includegraphics[scale=.38]{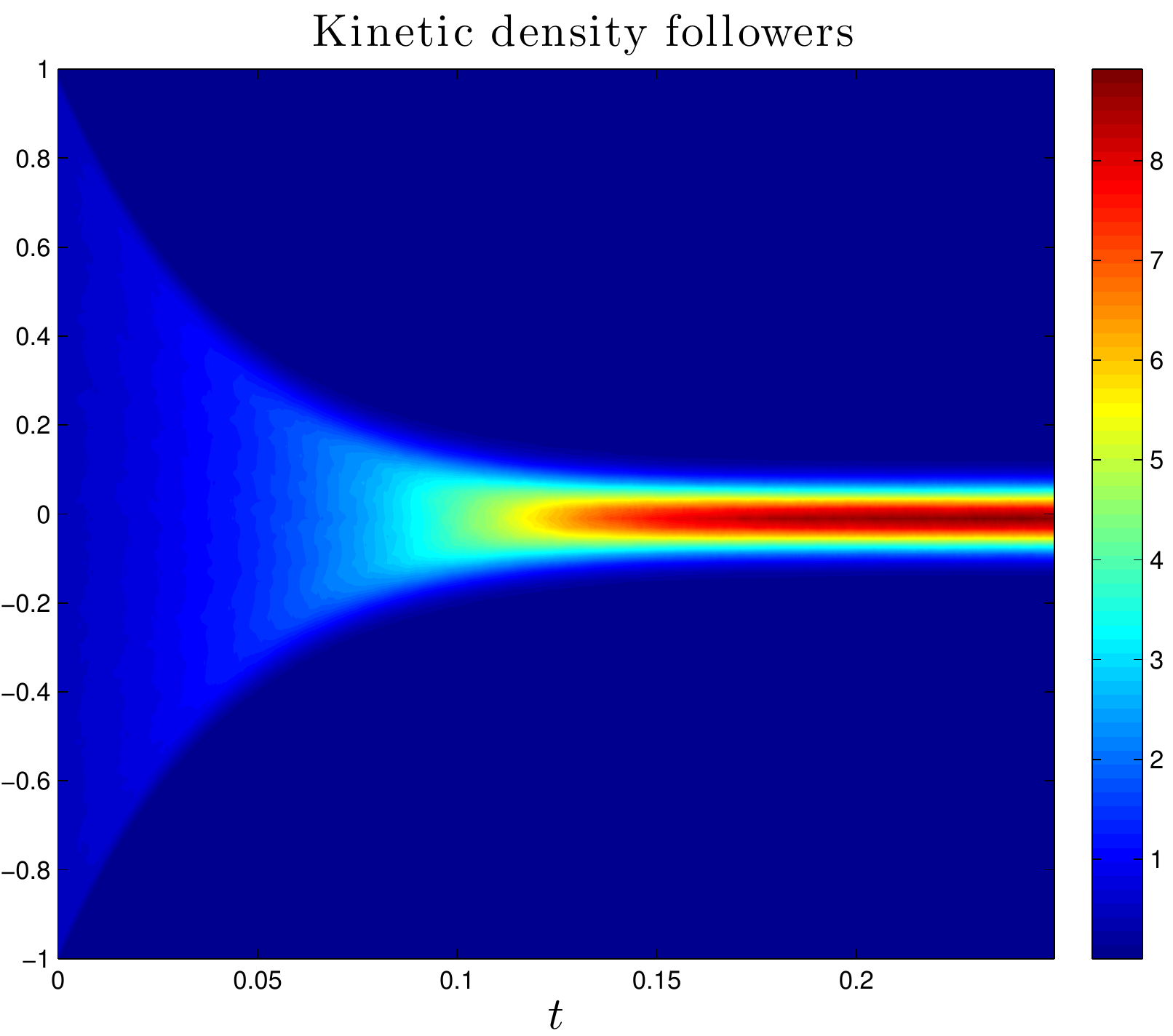}
\includegraphics[scale=.38]{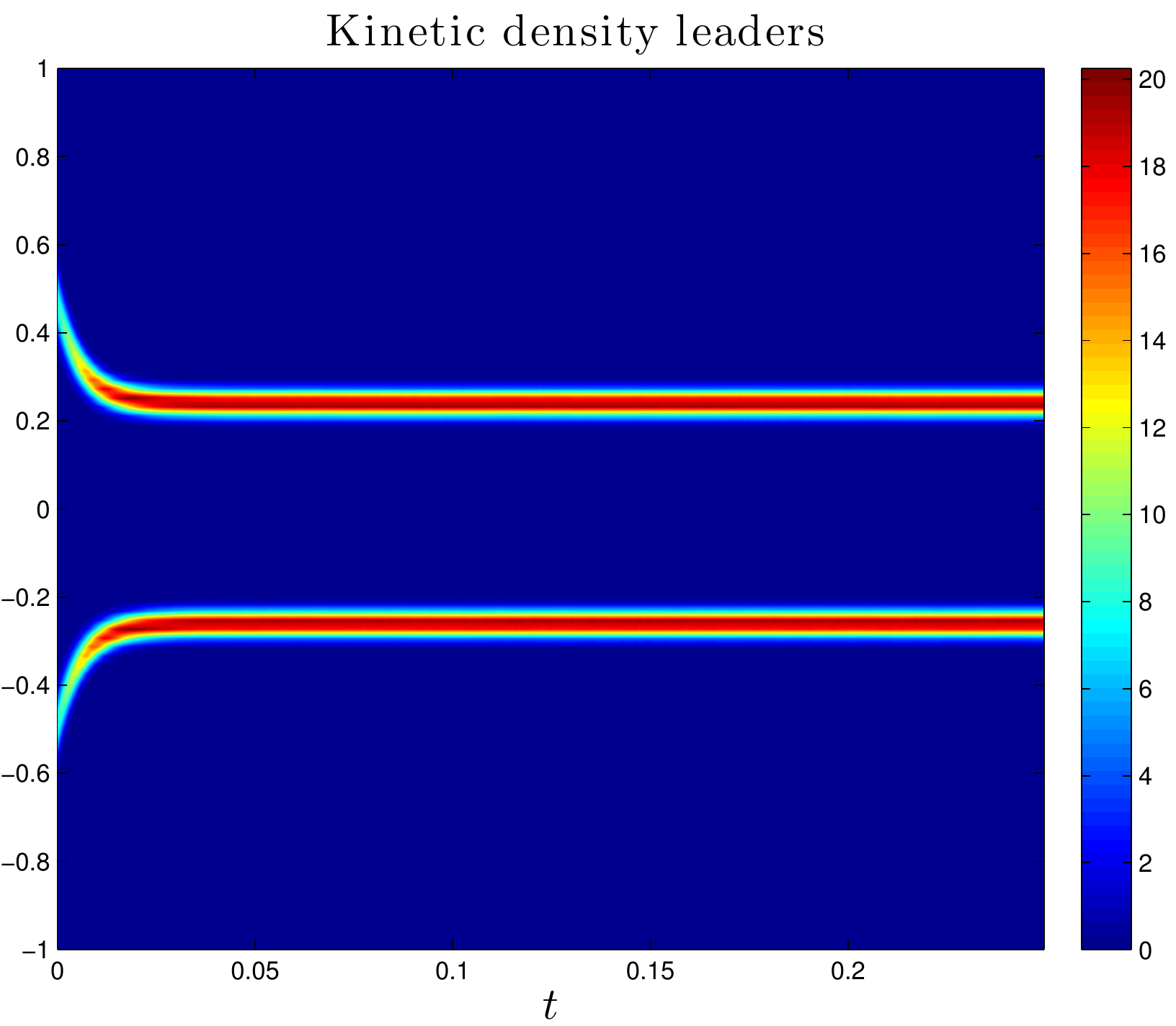}
\caption{Test 2: Kinetic densities evolution over the time interval $[0,0.25]$ reproducing a Hotelling-like model behaviour for two populations of leaders.}
\label{fig:hot2}
\end{figure}

\subsection*{Test 2. The case of multiple leaders populations}
Similarly, if more than one population of leaders occurs, each one with a different strategy, we can describe the evolution of the kinetic density of the system through a Boltzmann approach. Let $M> 0$ be the number of families of leaders, each of them described by the density $f_{L_p},p=1,...,M$ such that
\begin{equation}
\int_I f_{L_p}(\tilde{w})d\tilde{w}=\rho_{p}.
\end{equation}
If we suppose that an unique population of followers does exist, with density $f_F$, and that every follower interacts both with the others agents from the same population and with every leader of each $p$-th family, for a suitable test function $\varphi$ we obtain the following system of Boltzmann equations 
\begin{equation}\begin{cases} \vspace{2 mm}
\dfrac{d}{dt}\displaystyle\int_I \varphi(w)f_F(w,t)dw=\left(Q_F(f_F,f_F),\varphi\right)+\displaystyle\sum_{p=1}^{M}\left(Q_{FL}(f_{L_p},f_F),\varphi\right) &\\ 
\dfrac{d}{dt}\displaystyle\int_I\varphi(\tilde{w})f_{L_p}(\tilde{w},t)d\tilde{w}=(Q_L(f_{L_p},f_{L_p}),\varphi),&\quad p=1,\ldots,M.
\end{cases}\end{equation} 
 We assume that the leaders aim at minimising cost functionals of the type (\ref{integral_functional}) and therefore the differences consist in two factors: the target opinions $w_{d_p}$ and in the leaders’ attitude towards a radical ($\psi_p\approx 1$) or populistic strategy ($\mu_p\approx 1$). We therefore introduce the analogous rescaling \eqref{eq:scaling} and we define 
 \be
 \hat{c}_{FL_p}={c_{FL_p}}/{\rho_p},\quad \hat{c}_{L_p}={c_{L_p}}/{\rho_p},\quad p=1,\ldots,M.
 \ee  
 
In the numerical test we establish a link between our arguments and a Hotelling's type model \cite{Hotelling29}. The model describes how two shop owners, which sell the same product at the same price in the same street, must locate their shops in order to reach the maximum number of customers, uniformly distributed along the street (in other words, in order to maximize their profits). Paradoxically the model yields that the equilibrium, without changing prices, is reached if they get closer. In the cited original paper electoral dynamics are placed in this context and it can be regarded as the reason why political parties' programs are often perceived as similar. We consider the case of two populations of leaders, described by the densities $f_{L_1}$ and $f_{L_2}$, exercising different controls over a population of followers uniformly distributed within the interval $I=[-1,1]$. Initially the leaders are distributed as $f_{L_p}\sim N(w_{d_p},0.05)$, $p=1,2$. We can observe that the model leads to a centrist population of followers, whose opinion spreads in a range between leaders' mean opinions (see Figures \ref{fig:hot1} and \ref{fig:hot2}). 

\begin{figure}[t]
\centering
\includegraphics[scale=0.38]{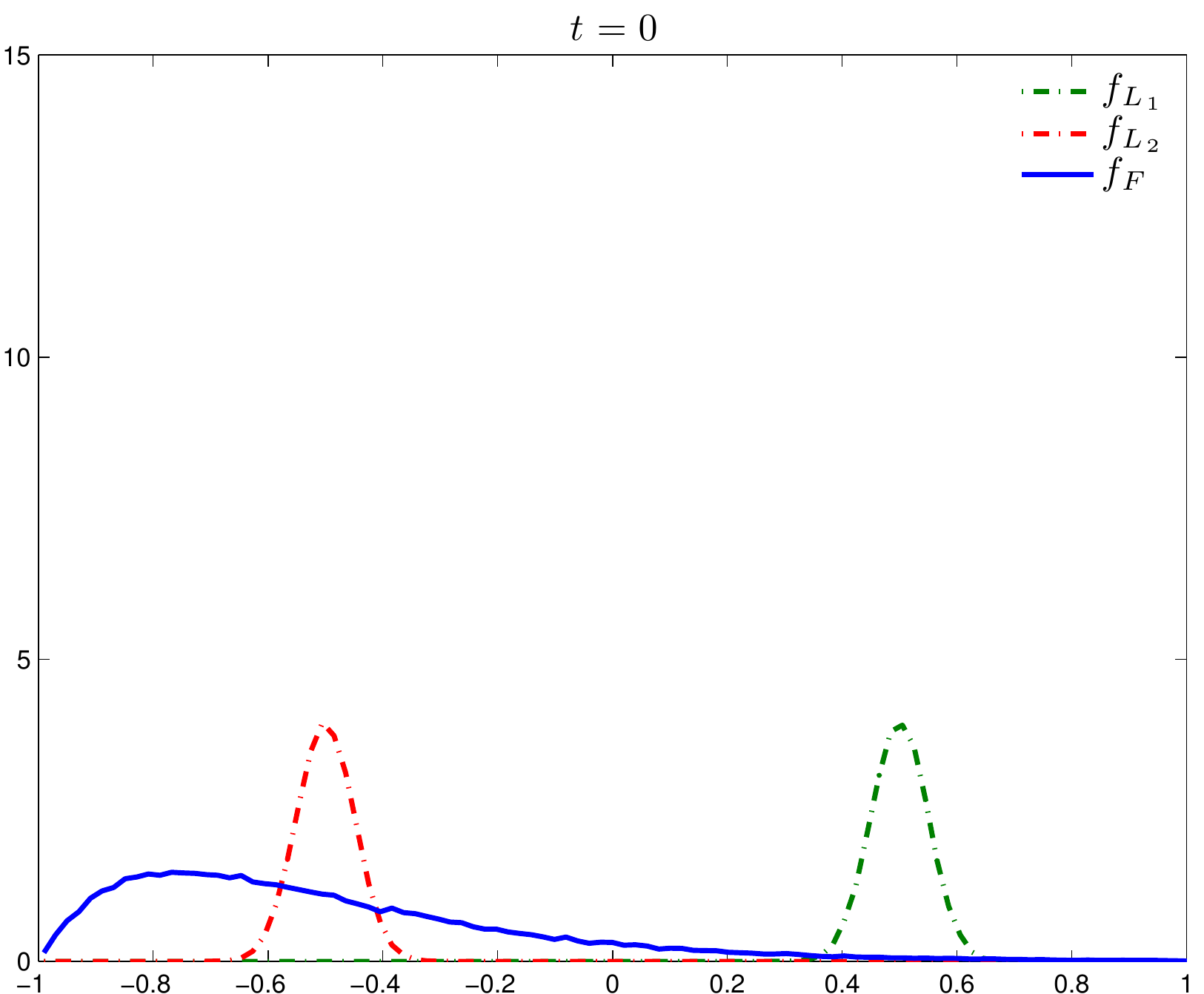}
\includegraphics[scale=0.38]{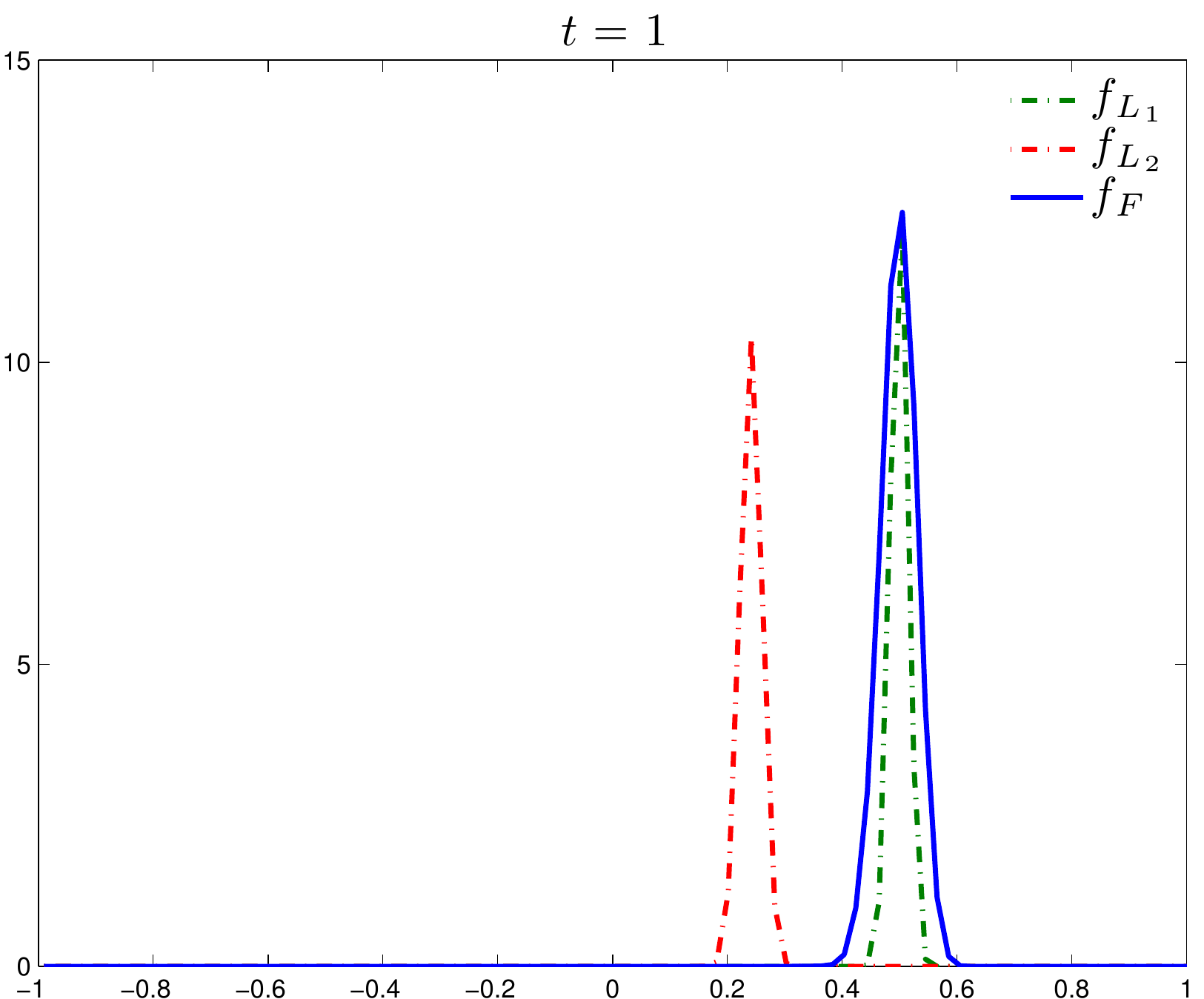}
\caption{Test 3: Kinetic densities at different times for  for a two populations of leaders model with time dependent strategies.}
\label{fig:2pop_time_dep}
\end{figure}

\begin{figure}[t]
\centering
\includegraphics[scale=0.38]{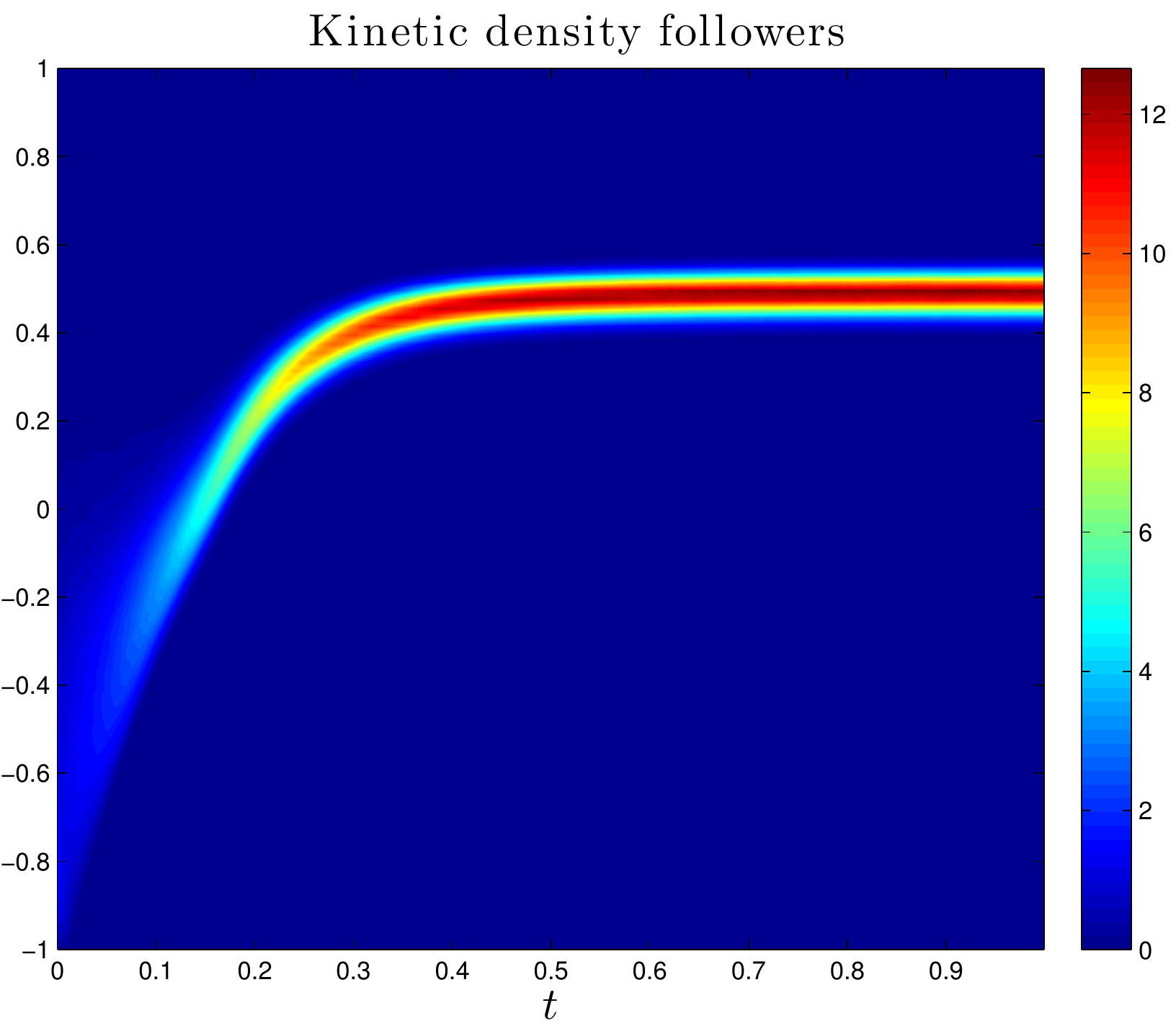}
\includegraphics[scale=0.38]{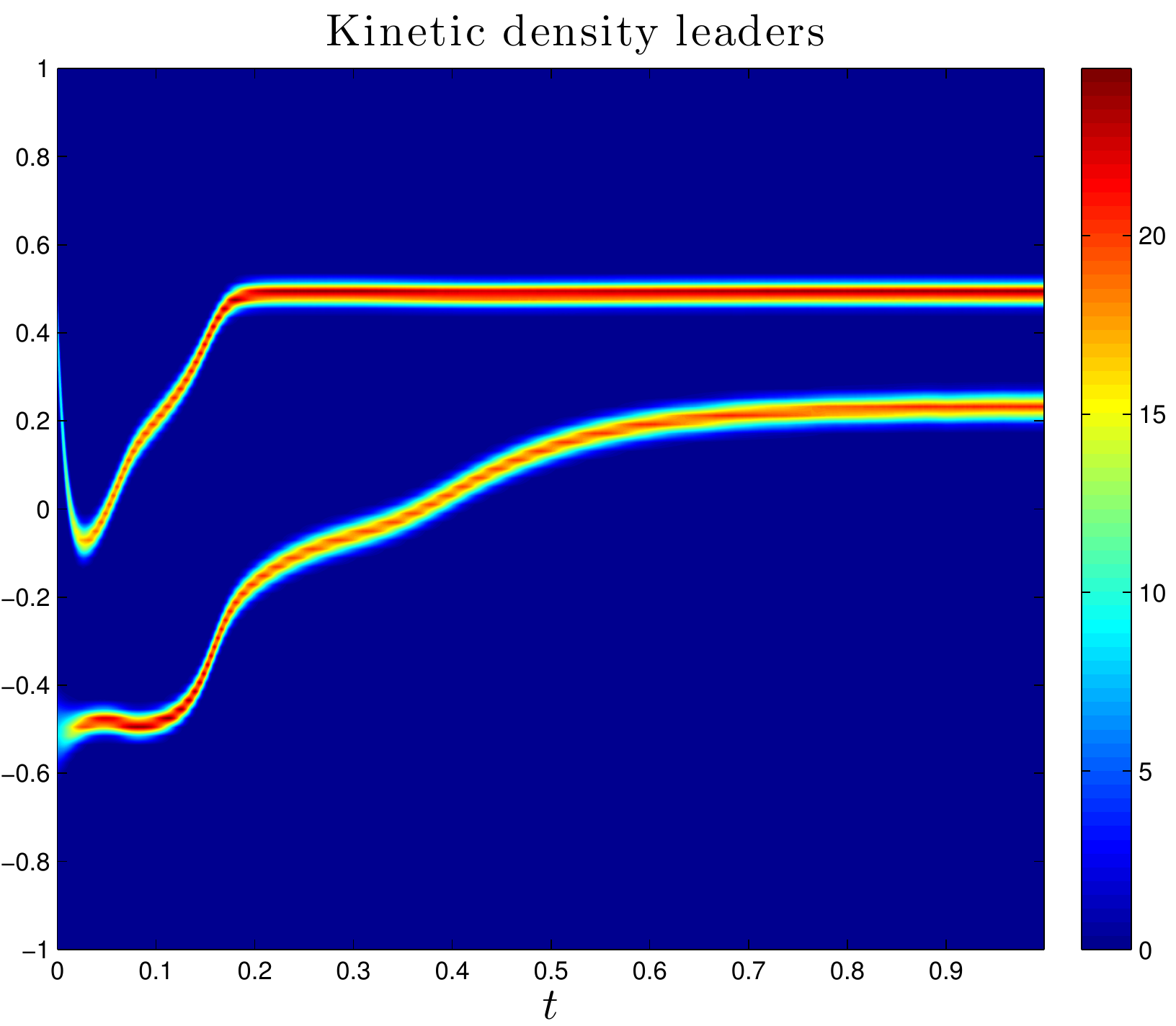}
\caption{Test 3: Kinetic density evolution over the time interval $[0,1]$ for a two populations of leaders model with time dependent strategies.}
\label{fig:kin_densities_timedep}
\end{figure}

\subsection*{Test 3. Two leaders populations with time-dependent strategies} Finally we introduce a multi-population model for opinion formation with time-dependent coefficients. 
This approach leads to the concept of adaptive strategy for every family of leaders $p=1,...,M$. The coefficients $\psi$ and $\mu$ which appear into the functional now evolve in time and are defined for every $t \in [0,T]$ as 
\begin{equation}\begin{split}\label{psi_timedep}
\psi_p(t)=&\dfrac{1}{2}{\int_{w_{d_p}-\delta}^{w_{d_p}+\delta}f_F(w)dw}+\dfrac{1}{2}\int_{{m_{L_p}}-\bar{\delta}}^{{m_{L_p}}+\bar{\delta}}f_F(w)dw \\
\mu_p(t)=&1-\psi_p(t)
\end{split}\end{equation}
where both $\delta,\bar{\delta} \in [0,1]$ are fixed and $m_{L_p}$ is the average opinion of the $p$-th leader. This choice of coefficients is equivalent to introduce a competition between the populations of leaders, where each leader try to adapt its populistic or radical attitude accordingly to the success of the strategy. Note also that the success of the strategy is based on the local perception of the followers.  

In the numerical experiments reported in Figure \ref{fig:2pop_time_dep} and Figure \ref{fig:kin_densities_timedep} we take into account two populations of leaders, initially normally distributed with mean values $w_{d_1}$ and $w_{d_2}$ and parameters $\delta=\bar{\delta}=0.5$, respectively, and a single population of followers, represented by a skewed distribution $f_F \sim \Gamma (2,\frac{1}{4})$ over the interval $[-1,1]$, where $\Gamma(\cdot,\cdot)$ is the Gamma distribution. Here the frequencies of interactions are assumed to be unbalanced since $\hat{c}_{FL_1}=0.1$ and $\hat{c}_{FL_2}=1$. In the test case we assume that the followers group has an initial natural inclination for a position represented by one leader but, thanks to communication strategies pursued by the minority leader, it is driven to different positions (see Figures \ref{fig:2pop_time_dep} and \ref{fig:kin_densities_timedep}). In a bipolar electoral context, an example of the described behaviour would be a better use of the media in a coalition with respect to the opponents. 

\section{Conclusions}
We introduced a Boltzmann type control for a hierarchical model of opinion formation where the leader behaviour is influenced both by the desire to achieve a prescribed opinion consensus and by the mean opinion of the followers. 
The main novelty of the method is that, thanks to an instantaneous binary control approximation, the control is explicitly incorporated in the resulting leader dynamic. The use of instantaneous control and the kinetic description permit to pass from an $O(N^2)$ dynamic, which must be solved forward-backward in time, to a much simpler forward $O(N)$ stochastic simulation. This is of paramount importance in view of possible applications of this kind of constrained opinion modeling. In the so-called quasi invariant opinion limit the corresponding Fokker-Planck descriptions have been derived and explicit expressions of their steady states computed. Several numerical examples illustrate the robustness of the controlled dynamics using various leaders strategies even in presence of different groups of competing leaders.

\bigskip
\noindent
{\bf Acknowledgements.} This work has been supported
by PRIN-MIUR grant {\em Advanced numerical methods for
kinetic equations and balance laws with source terms} and by FAR 2013 project {\em Multiscale PDEs and emerging collective behaviours}. G. Albi acknowledges the support of the ERC-Starting Grant HDSPCONTR {\em High-Dimensional Sparse Optimal Control}. 


\end{document}